\documentclass[final]{siamart0516}

\usepackage{lipsum}
\usepackage{amsfonts,amssymb}
\usepackage{graphicx}
\usepackage{epstopdf}
\usepackage{algorithmic}
\usepackage{tikz}
\usetikzlibrary{shapes,arrows}

\input{mysymbol.sty}

\ifpdf
  \DeclareGraphicsExtensions{.eps,.pdf,.png,.jpg}
\else
  \DeclareGraphicsExtensions{.eps}
\fi

\newcommand{\TheTitle}{{IQN: An Incremental Quasi-Newton Method \\ with Local Superlinear Convergence Rate}} 
\newcommand{\TheAuthors}{Aryan Mokhtari, Mark Eisen, and Alejandro Ribeiro}

\headers{IQN: An Incremental Quasi-Newton Method}{\TheAuthors}

\title{{\TheTitle}\thanks{{This work was supported by ONR N00014-12-1-0997.}}}

\author{ Aryan Mokhtari\thanks{Department of Electrical and Systems Engineering, University of Pennsylvania, Philadelphia, PA
    (\email{aryanm@seas.upenn.edu},\email{maeisen@seas.upenn.edu},\email{aribeiro@seas.upenn.edu})}
  \and
 Mark Eisen\footnotemark[2]
  \and
  Alejandro Ribeiro\footnotemark[2]
}

\usepackage{amsopn}

\newtheorem{assumption}{\hspace{0pt}\bf Assumption}
\newtheorem{remark}{\hspace{0pt}\it Remark}

\begin{document}

\maketitle

\begin{abstract}
The problem of minimizing an objective that can be written as the sum of a set of $n$ smooth and strongly convex functions is challenging because the cost of evaluating the function and its derivatives is proportional to the number of elements in the sum. The Incremental Quasi-Newton (IQN) method proposed here belongs to the family of stochastic and incremental methods that have a cost per iteration independent of $n$. IQN iterations are a stochastic version of BFGS iterations that use memory to reduce the variance of stochastic approximations. The method is shown to exhibit local superlinear convergence. The convergence properties of IQN bridge a gap between deterministic and stochastic quasi-Newton methods. Deterministic quasi-Newton methods exploit the possibility of approximating the Newton step using objective gradient differences. They are appealing because they have a smaller computational cost per iteration relative to Newton's method and achieve a superlinear convergence rate under customary regularity assumptions. Stochastic quasi-Newton methods utilize stochastic gradient differences in lieu of actual gradient differences. This makes their computational cost per iteration independent of the number of objective functions $n$. However, existing stochastic quasi-Newton methods have sublinear or linear convergence at best. IQN is the first stochastic quasi-Newton method proven to converge superlinearly in a local neighborhood of the optimal solution. IQN differs from state-of-the-art incremental quasi-Newton methods in three aspects: (i) The use of aggregated information of variables, gradients, and quasi-Newton Hessian approximation matrices to reduce the noise of gradient and Hessian approximations. (ii) The approximation of each individual function by its Taylor's expansion in which the linear and quadratic terms are evaluated with respect to the same iterate. (iii) The use of a cyclic scheme to update the functions in lieu of a random selection routine. We use these fundamental properties of IQN to establish its local superlinear convergence rate. The presented numerical experiments match our theoretical results and justify the advantage of IQN relative to other incremental methods.
\end{abstract}


\begin{AMS}
{  90C06, 90C25, 90C30, 90C52}
\end{AMS}

\section{Introduction} \label{sec_intro}

We study large scale optimization problems with objective functions expressed as the sum of a set of components which arise often in application domains such as machine learning~\cite{BottouCun, Bottou, SS,cevher2014convex}, control \cite{Bullo2009, Cao2013-TII, LopesEtal8}, and wireless communications \cite{Schizas2008-1, Ribeiro10, Ribeiro12}. Formally, we consider a variable $\bbx\in\reals^p$ and a function $f$ which is defined as the average of $n$ smooth and strongly convex functions labelled $f_i: \mathbb{R}^p \rightarrow \mathbb{R}$ for $i=1,\hdots, n$. We refer to individual functions $f_i$ as sample functions and to the total number of functions $n$ as the sample size. Our goal is to find the optimal argument $\bbx^*$ that solves the strongly convex program
\begin{equation}\label{eq_orig_problem}
   \bbx^* := \argmin_{\bbx \in \reals^p}  f(\bbx) 
          := \argmin_{\bbx \in \reals^p} \frac{1}{n} \sum_{i=1}^n f_i(\bbx).
\end{equation}
We restrict attention to cases where the component functions $f_i$ are strongly convex and their gradients are Lipschitz continuous. We further focus in problems where $n$ is large enough so as to warrant application of stochastic or iterative methods. Our goal is to propose an iterative quasi-Newton method to solve \eqref{eq_orig_problem} which is shown to exhibit a local superlinear convergence rate. This is achieved while performing local iterations with a cost of order $\mathcal{O}(p^2)$ independent of the number of samples $n$.

Setting temporarily aside the complications related to the number of component functions, the minimization of $f$ in \eqref{eq_orig_problem} can be carried out using iterative descent algorithms. A simple solution is to use gradient descent (GD) which iteratively descends along gradient directions $\nabla f(\bbx) = (1/n) \sum_{i=1}^n \nabla f_i(\bbx)$. GD incurs a per iteration computational cost of order $\mathcal{O}(np)$ and is known to converge at a linear rate towards $\bbx^*$ under the hypotheses we have placed on $f$. Whether the linear convergence rate of GD is acceptable depends on the desired accuracy and on the condition number of $f$ which, when large, can make the convergence constant close to one. As one or both of these properties often limit the applicability of GD, classical alternatives to improve convergence rates have been developed. Newton's method adapts to the curvature of the objective by computing Hessian inverses and converges at a quadratic rate in a local neighborhood of the optimal argument irrespective of the problem's condition number. To achieve this quadratic convergence rate, we must evaluate and invert Hessians resulting in a per iteration cost of order $\mathcal{O}(np^2+p^3)$. Quasi-Newton methods build on the idea of approximating the Newton step using first-order information of the objective function and exhibit local superlinear convergence \cite{Broyden, Powell, Dennis}. An important feature of quasi-Newton methods is that they have a per iteration cost of order $\mathcal{O}(np+p^2)$, where the term $\mathcal{O}(np)$ corresponds to the cost of gradient computation and the cost $\mathcal{O}(p^2)$ indicates the computational complexity of updating the approximate Hessian inverse matrix. 

The combination of a local superlinear convergence rate and the smaller computational cost per iteration relative to Newton -- a reduction by a factor of $p$ operations per iteration -- make quasi-Newton methods an appealing choice. In the context of optimization problems having the form in \eqref{eq_orig_problem}, quasi-Newton methods also have the advantage that curvature is estimated using gradient evaluations. To see why this is meaningful we must recall that the customary approach to avoid the $\mathcal{O}(np)$ computational cost of GD iterations is to replace gradients $\nabla f(\bbx)$ by their stochastic approximations $\nabla f_i(\bbx)$, which can be evaluated with a cost of order $\mathcal{O}(p)$. One can then think of using stochastic versions of these gradients to develop stochastic quasi-Newton methods with per iterations cost of order $\mathcal{O}(p+p^2)$. This idea was demonstrated to be feasible in \cite{SchraudolphYG07} which introduces a stochastic (online) version of the BFGS quasi-Newton method as well as a stochastic version of its limited memory variant. Although \cite{SchraudolphYG07} provides numerical experiments illustrating significant improvements in convergence times relative to stochastic (S) GD, theoretical guarantees are not established. 

The issue of proving convergence of stochastic quasi-Newton methods is tackled in \cite{mokhtari2014res} and \cite{mokhtari2015global}. In \cite{mokhtari2014res} the authors show that stochastic BFGS may not be convergent because the Hessian approximation matrices can become close to singular. A regularized stochastic BFGS (RES) method is proposed by changing the proximity condition of BFGS to ensure that the eigenvalues of the Hessian inverse approximation are uniformly bounded. Enforcing this property yields a provably convergent algorithm. In \cite{mokhtari2015global} the authors show that the limited memory version of stochastic (online) BFGS proposed in \cite{SchraudolphYG07} is almost surely convergent and has a sublinear convergence rate in expectation. This is achieved without using regularizations. An alternative provably convergent stochastic quasi-Newton method is proposed in \cite{byrd2016stochastic}. This method differs from those in \cite{SchraudolphYG07, mokhtari2014res, mokhtari2015global} in that it collects (stochastic) second order information to estimate the objective's curvature. This is in contrast to estimating curvature using the difference of two consecutive stochastic gradients.

Although the methods in \cite{SchraudolphYG07, mokhtari2014res, mokhtari2015global, byrd2016stochastic} are successful in expanding the application of quasi-Newton methods to stochastic settings, their convergence rate is sublinear. This is not better than the convergence rate of SGD and, as is also the case in SGD, is a consequence of the stochastic approximation noise which necessitates the use of diminishing stepsizes. The stochastic quasi-Newton methods in \cite{lucchi2015variance, MoritzNJ16} resolve this issue by using the variance reduction technique proposed in \cite{johnson2013accelerating}. The fundamental idea of the work in \cite{johnson2013accelerating} is to reduce the noise of the stochastic gradient approximation by computing the exact gradient in an outer loop to use it in an inner loop for gradient approximation. The methods in  \cite{lucchi2015variance, MoritzNJ16}, which incorporate the variance reduction scheme presented in \cite{johnson2013accelerating} into the update of quasi-Newton methods, are successful in achieving a linear convergence rate. 

At this point, we must remark on an interesting mismatch. The convergence rate of SGD is sublinear, and the convergence rate of deterministic GD is linear. The use of variance reduction techniques in SGD recovers the linear convergence rate of GD, \cite{johnson2013accelerating}. On the other hand, the convergence rate of stochastic quasi-Newton methods is sublinear, and the convergence rate of deterministic quasi-Newton methods is superlinear. The use of variance reduction in stochastic quasi-Newton methods achieves linear convergence but does not recover a superlinear rate. Hence, a fundamental question remains unanswered: Is it possible to design an incremental quasi-Newton method that recovers the superlinear convergence rate of deterministic quasi-Newton algorithms? In this paper, we show that the answer to this open problem is positive by proposing an incremental quasi-Newton method (IQN) with a local superlinear convergence rate. This is the first quasi-Newton method to achieve superlinear convergence while having a per iteration cost independent of the number of functions $n$ -- the cost per iteration is of order $\mathcal{O}(p^2)$.

There are three major differences between the IQN method and state-of-the-art incremental (stochastic) quasi-Newton methods that lead to the former's superlinear convergence rate. First, the proposed IQN method uses the aggregated information of variables, gradients, and Hessian approximation matrices to reduce the noise of approximation for both gradients and Hessian approximation matrices. This is different to the variance-reduced stochastic quasi-Newton methods in \cite{lucchi2015variance,MoritzNJ16} that attempt to reduce only the noise of gradient approximations. Second, in IQN the index of the updated function is chosen in a cyclic fashion, rather than the random selection scheme used in the incremental methods in \cite{SchraudolphYG07, mokhtari2014res, mokhtari2015global, byrd2016stochastic}. The cyclic routine in IQN allows to bound the error at each iteration as a function of the errors of the last $n$ iterates, something that is not possible when using a random scheme. To explain the third and most important difference we point out that the form of quasi-Newton updates is the solution of a local second order Taylor approximation of the objective. It is possible to understand stochastic quasi-Newton methods as an analogous approximation of individual sample functions. However, it turns out that the state-of-the-art stochastic quasi-methods evaluate the linear and quadratic terms of the Taylor's expansion at different points yielding and inconsistent approximation (Remark \ref{eq_second_approx}). The IQN method utilizes a consistent Taylor series which yields a more involved update which we nonetheless show can be implemented with the same computational cost. These three properties together lead to an incremental quasi-Newton method with a local superlinear convergence rate. 

We start the paper by recapping the BFGS quasi-Newton method and the Dennis-Mor\'e condition which is sufficient and necessary to prove superlinear convergence rate of the BFGS method (Section~\ref{sec_problem_formulation1}). Then, we present the proposed Incremental Quasi-Newton method (IQN) as an incremental aggregated version of the traditional BFGS method (Section~\ref{sec_inc_bfgs}). We first explain the difference between the Taylor's expansion used in IQN and state-of-the-art incremental (stochastic) quasi-Newton methods. Further, we explain the mechanism for aggregation of the functions information and the scheme for updating the stored information. Moreover, we present an efficient implementation of the proposed IQN method with a computational complexity of the order $\mathcal{O}(p^2)$ (Section~\ref{sec_iqn_implementation}). The convergence analysis of the IQN method is then presented (Section~\ref{sec_conv_analysis}). We use the classic analysis of quasi-Newton methods to show that in a local neighborhood of the optimal solution the sequence of variables converges to the optimal argument $\bbx^*$ linearly after each pass over the set of functions (Lemma~3). We use this result to show that for each component function $f_i$ the Dennis-Mor\'e condition holds (Proposition \ref{prop_DM_condition}). However, this condition is not sufficient to prove superlinear convergence of the sequence of errors $\|\bbx^t-\bbx^*\|$, since it does not guarantee the Dennis-Mor\'e condition for the global objective $f$. To overcome this issue we introduce a novel convergence analysis approach which exploits the local linear convergence of IQN to present a more general version of the Dennis-Mor\'e condition for each component function $f_i$ (Lemma \ref{lemma_incr_bfgs_lim_final}). We exploit this result to establish superlinear convergence of the iterates generated by IQN (Theorem \ref{sup_linear_thm}). In Section~\ref{sec_num_results}, we present numerical simulation results, comparing the performance of IQN to that of first-order incremental and stochastic methods. We test the performance on a set of large-scale regression problems and observe strong numerical gain in total computation time relative to existing methods.

\subsection{Notation}Vectors are written as lowercase $\bbx\in\reals^p$ and matrices as uppercase $\bbA\in\reals^{p\times p}$. We use $\|\bbx\|$ and  $\|\bbA\|$ to denote the Euclidean norm of vector $\bbx$ and matrix $\bbA$, respectively. Given a positive definite matrix $\bbM$, the weighted matrix norm $\|\bbA\|_\bbM$ is defined as $\|\bbA\|_\bbM:= \|\bbM\bbA\bbM\|_\bbF$, where $\|.\|_\bbF$ is the Frobenius norm. Given a function $f$ its gradient and Hessian at point $\bbx$ are denoted as $\nabla f(\bbx)$ and $\nabla^2 f(\bbx)$, respectively. 


\section{BFGS Quasi-Newton Method}\label{sec_problem_formulation1}

Consider the problem in \eqref{eq_orig_problem} for relatively large $n$. In a conventional optimization setting, this can be solved using a quasi-Newton method that iteratively updates a variable $\bbx^t$ for $t=0,1,\hdots$ based on the general recursive expression
\begin{align} \label{eq_descent_update}
\bbx^{t+1} = \bbx^t - \eta^t({\bbB^t})^{-1}\nabla f(\bbx^t),
\end{align}
where $\eta^t$ is a scalar stepsize and ${\bbB^t}$ is a positive definite matrix that approximates the exact Hessian of the objective function $\nabla^{2} f(\bbx^t)$. The stepsize $\eta^t$ is evaluated based on a line search routine for the global convergence of quasi-Newton methods. Our focus in this paper, however, is on the local convergence of quasi-Newton methods, which requires the unit stepsize $\eta^t=1$. Therefore, throughout the paper we assume that the variable $\bbx^t$ is close to the optimal solution $\bbx^*$ -- we will formalize the notion of being close to the optimal solution -- and the stepsize is $\eta^t=1$.

 The goal of quasi-Newton methods is to compute the Hessian approximation matrix  ${\bbB^t}$ and its inverse ${(\bbB^t)}^{-1}$ by using only the first-order information, i.e., gradients, of the objective. Their use is widespread due to the many applications in which the Hessian information required in Newton's method is either unavailable or computationally intensive. There are various approaches to approximate the Hessian, but the common feature among quasi-Newton methods is that the Hessian approximation must satisfy the secant condition. To be more precise, consider $\bbs^t $ and $\bby^t$ as the variable and gradient variations, explicitly defined as
%
%
%
%
\begin{align}\label{eq_bfgs_vars}
\bbs^t := \bbx^{t+1} - \bbx^t, \qquad \bby^t := \nabla f(\bbx^{t+1}) - \nabla f(\bbx^t).
\end{align}
Then, given the variable variation $\bbs^t$ and gradient variation $\bby^t$, the Hessian approximation matrix in all quasi-Newton methods must satisfy the secant condition
\begin{equation}\label{secant_condition}
\bbB^{t+1}\bbs^t=\bby^t.
\end{equation}
This condition is fundamental in quasi-Newton methods because the exact Hessian $\nabla^{2} f(\bbx^t)$ satisfies this equality when the iterates $\bbx^{t+1}$ and $\bbx^t$ are close to each other. If we consider the matrix $\bbB^{t+1}$ as the unknown matrix, the system of equations in \eqref{secant_condition} does not have a unique solution. Different quasi-Newton methods enforce different conditions on the matrix $\bbB^{t+1}$ to come up with a unique update. This extra condition is typically a proximity condition that ensures that $\bbB^{t+1}$ is close to the previous Hessian approximation matrix $\bbB^t$ \cite{Broyden,Powell,Dennis}. In particular, the Broyden-Fletcher-Goldfarb-Shanno (BFGS) method defines the update of Hessian approximation matrix as
\begin{align}\label{eq_bfgs_update}
\bbB^{t+1} = \bbB^t + \frac{\bby^t {\bby^t}^T}{{\bby^t}^T \bbs^t} - \frac{ \bbB^t \bbs^t {\bbs^t}^T \bbB^t}{{\bbs^t}^T \bbB^t \bbs^t}.
\end{align}
The BFGS method is popular not only for its strong numerical performance relative to the gradient descent method, but also because it is shown to exhibit a superlinear convergence rate \cite{Broyden}, thereby providing a theoretical guarantee of superior performance. In fact, it can be shown that, the BFGS update satisfies the condition
\begin{align}\label{eq_dennis_more}
\lim_{t \rightarrow \infty} \frac{ \| (\bbB^t - \nabla^2 f(\bbx^*)) \bbs^t \|}{\|\bbs^t\|} = 0,
\end{align}
known as the Dennis-Mor\'e condition, which is both necessary and sufficient for superlinear convergence \cite{Dennis}. This result solidifies quasi-Newton methods as a strong alternative to first order methods when exact second-order information is unavailable. However, implementation of the BFGS method is not feasible when the number of functions $n$ is large, due to its high computational complexity on the order $\mathcal{O}(np+p^2)$. In the following section, we propose a novel incremental BFGS method that has the computational complexity of $\mathcal{O}(p^2)$ per iteration and converges at a superlinear  rate.

%

\section{IQN: Incremental aggregated BFGS}\label{sec_inc_bfgs}

We propose an incremental aggregated BFGS algorithm, which we call the Incremental Quasi-Newton (IQN) method. The IQN method is incremental in that, at each iteration, only the information associated with a single function $f_i$ is updated. The particular function is chosen by cyclicly iterating through the $n$ functions. The IQN method is aggregated in that the aggregate of the most recently observed information of all functions $f_1,\dots,f_n$ is used to compute the updated variable $\bbx^{t+1}$.

\begin{figure} \centering
\input{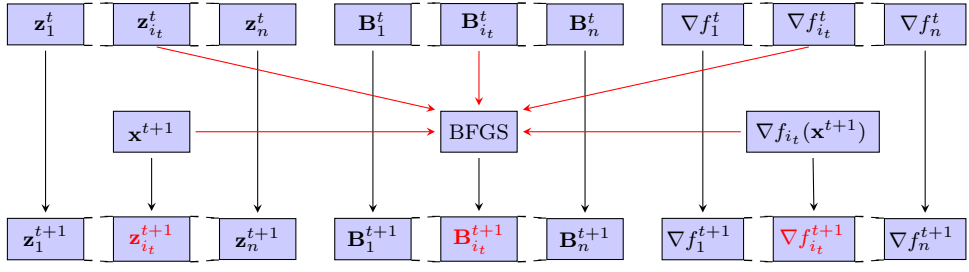}
\caption{The updating scheme for variables, gradients, and Hessian approximation matrices of function $f_{i_t}$ at step $t$. The red arrows indicate the terms used in the update of $\bbB_{i_t}^{t+1}$ using the BFGS update in \eqref{eq_bfgs_update_local}. The black arrows show the updates of all variables and gradients. The terms $\bbz_{i_t}^{t+1}$ and $\nabla f_{i_t}^{t+1}$ are updated as $\bbx^{t+1}$ and $\nabla f_{i_t}(\bbx^{t+1})$, respectively. All others $\bbz_{j}^{t+1}$ and $\nabla f_{j}^{t+1}$ are set as $\bbz_{j}^{t}$ and $\nabla f_{j}^{t}$, respectively.  }
\label{fig_sag} \end{figure}

 In the proposed method, we consider $\bbz_1^t, \hdots, \bbz_n^t$ as the copies of the variable $\bbx$ at time $t$ associated with the functions $f_1,\dots,f_n$, respectively. Likewise, define $\nabla f_i(\bbz_i^t)$ as the gradient corresponding to the $i$-th function. Further, consider $\bbB_i^t$ as a positive definite matrix which approximates the $i$-th component Hessian $\nabla^2 f_i (\bbx^t)$. We refer to $\bbz_i^t$, $\nabla f_i(\bbz_i^t)$, and $\bbB_i^t$ as the information corresponding to the $i$-th function $f_i$ at step $t$. Note that the functions' information is stored in a shared memory as shown in Fig.~\ref{fig_sag}. To introduce the IQN method, we first explain the mechanism for computing the updated variable $\bbx^{t+1}$ using the stored information $\{\bbz_i^t,\nabla f_i(\bbz_i^t),\bbB_i^t\}_{i=1}^n$. Then, we elaborate on the scheme for updating the information of the functions.


To derive the full variable update, consider the second order approximation of the objective function $f_i(\bbx)$ centered around its current iterate $\bbz_i^t$,
\begin{align}\label{eq_second_approx}
f_i(\bbx) \approx f_i (\bbz_i^t) + \nabla f_i (\bbz_i^t)^T (\bbx - \bbz_i^t) 
 + \frac{1}{2} (\bbx- \bbz_i^t)^T \nabla^2 f_i (\bbz_i^t)(\bbx- \bbz_i^t).
\end{align}
As in traditional quasi-Newton methods, we replace the $i$-th Hessian $\nabla^2 f_i (\bbz_i^t)$ by $\bbB_i^t$. Using the approximation matrices in place of Hessians, the complete (aggregate) function $f(\bbx)$ can be approximated with
\begin{align}\label{eq_second_approx_2}
f(\bbx) \approx \frac{1}{n} \sum_{i=1}^n \left[ f_i (\bbz_i^t) +  \nabla f_i (\bbz_i^t)^T (\bbx - \bbz_i^t)  
+ \frac{1}{2} (\bbx - \bbz_i^t)^T \bbB_i^t (\bbx- \bbz_i^t)\right].
\end{align}
Note that the right hand side of \eqref{eq_second_approx_2} is a quadratic approximation of the function $f$ based on the available information at step $t$. Hence, the updated iterate $\bbx^{t+1}$ can be defined as the minimizer of the quadratic program in \eqref{eq_second_approx_2}, explicitly given by
\begin{align}\label{eq_incr_bfgs_update_222}
\bbx^{t+1} = \left( \frac{1}{n} \sum_{i=1}^n \bbB_i^t \right)^{-1} \left[ \frac{1}{n} \sum_{i=1}^n \bbB_i^t \bbz_i^t - \frac{1}{n} \sum_{i=1}^n \nabla f_i (\bbz_i^t) \right].
\end{align}

First note that the update in \eqref{eq_incr_bfgs_update_222} shows that the updated variable $\bbx^{t+1}$ is a function of the stored information of all functions $f_1,\dots,f_n$. Furthermore, we use the aggregated information of variables, gradients, and the quasi-Newton Hessian approximations to evaluate the updated variable. This is done to vanish the noise in approximating both gradients and Hessians as the sequence approaches the optimal argument. 

\begin{remark}
Given the BFGS Hessian approximation matrices $\{\bbB_i^t\}_{i=1}^n$ and gradients $\{\nabla f_i (\bbz_i^t)\}_{i=1}^n$, one may consider an update more akin to traditional descent-based methods, i.e.
\begin{align}\label{eq_incr_bfgs_update_alternative}
\bbx^{t+1} = \bbx^t - \left( \frac{1}{n} \sum_{i=1}^n \bbB_i^t \right)^{-1} \frac{1}{n} \sum_{i=1}^n \nabla f_i (\bbz_i^t) .
\end{align}
To evaluate the advantage of the proposed update for IQN in \eqref{eq_incr_bfgs_update_222} relative to the update in \eqref{eq_incr_bfgs_update_alternative}, we proceed to study the Taylor's expansion that leads to the update in \eqref{eq_incr_bfgs_update_alternative}. It can be shown that the update in \eqref{eq_incr_bfgs_update_alternative} is the outcome of the following approximation 
\begin{align}\label{eq_second_approx__alterantive}
f(\bbx) \approx \frac{1}{n} \sum_{i=1}^n \left[ f_i (\bbz_i^t) +  \nabla f_i (\bbz_i^t)^T (\bbx - \bbz_i^t)  
+ \frac{1}{2} (\bbx - \bbx^t)^T \bbB_i^t (\bbx- \bbx^t)\right].
\end{align}
Observe that the linear term in \eqref{eq_second_approx__alterantive} is centered at $\bbz_i^t$, while the quadratic term is approximated near the iterate $\bbx^t$. This inconsistency in the Taylor's expansion of each function $f_i$ leads to an inaccurate second-order approximation, and subsequently a slower incremental quasi-Newton method.
\end{remark}


 

Thus far we have discussed the procedure to compute the updated variable $\bbx^{t+1}$ given the local iterates, gradients, and Hessian approximations at time $t$. Now, it remains to show how we update the local information of functions  $f_1,\dots,f_n$ using the variable $\bbx^{t+1}$. In each iteration of the IQN method, we update the local information of only a single function, chosen in a cyclic manner. Defining $i_t$ to be the index of the function selected at time $t$, we update the local variables $\bbz_{i_t}^{t+1}$, $ \nabla f_{i_t}(\bbz_i^{t+1})$, and $\bbB_i^{t+1}$ using the updated variable $\bbx^{t+1}$ while all other local variables remain unchanged. In particular, the variables $\bbz_{i}$ are updated as
\begin{align} \label{eq_local_variable_update_222}
\bbz_{i_t}^{t+1} = \bbx^{t+1}, \qquad \bbz_i^{t+1} = \bbz_i^{t} \quad \text{ for all } i \neq i_t.
\end{align}
Observe in the update in \eqref{eq_local_variable_update_222} that the variable associated with the function $f_{i_t}$ is set to be the updated variable $\bbx^{t+1}$ while the other iterates are simply kept as their previous value. Likewise, we update the table of gradients accordingly with the gradient of $f_{i_t}$ evaluated at the new variable $\bbx^{t+1}$. The rest of gradients stored in the memory will stay unchanged, i.e., 
\begin{align} \label{eq_local_gradient_update_333}
\nabla f_{i_t}(\bbz_i^{t+1})=\nabla f_{i_t}(\bbx^{t+1}), \qquad \nabla f_{i}(\bbz_i^{t+1})= \nabla f_{i}(\bbz_i^{t}) \quad \text{ for all } i \neq i_t.
\end{align}

To update the curvature information, it would be ideal to compute the Hessian $\nabla^2 f_{i_t}(\bbx^{t+1})$ and update the curvature information following the schemes for variables in \eqref{eq_local_variable_update_222} and gradients in \eqref{eq_local_gradient_update_333}. However, our focus is on the applications that the computation of the Hessian is either impossible or computationally expensive. Hence, to the update curvature approximation matrix $\bbB_{i_t}^t$ corresponding to the function $f_{i_t}$, we use the steps of BFGS in \eqref{eq_bfgs_update}. To do so, we define variable and gradient variations associated with each individual function $f_i$ as 
\begin{align} \label{eq_bfgs_vars_local_2}
\bbs_i^t := \bbz_i^{t+1} - \bbz_i^t, \qquad \bby_i^t := \nabla f_i(\bbz_i^{t+1}) - \nabla f_i(\bbz_i^t),
\end{align}
respectively. The Hessian approximation $ \bbB_{i_t}^t $ corresponding to the function $f_{i_t}$ can be computed using the update of BFGS as
\begin{align} \label{eq_bfgs_update_local}
\bbB_i^{t+1} = \bbB_i^t + \frac{\bby_i^t \bby_i^{tT}}{\bby_i^{tT} \bbs_i^t} - \frac{ \bbB_i^t \bbs_i^t \bbs_i^{tT} \bbB_i^t}{\bbs_i^{tT} \bbB_i^t \bbs_i^t}, \quad  \text{ for}\quad i=i_t.
\end{align}
Again, the Hessian approximation matrices for all other functions remain unchanged, i.e., $\bbB_i^{t+1} = \bbB_i^t$ for $ i \neq i_t $. The system of updates in \eqref{eq_local_variable_update_222}-\eqref{eq_bfgs_update_local} explains the mechanism of updating the information of the function $f_{i_t}$ at step $t$. Notice that to update the Hessian approximation matrix for the $i_t$-th function there is no need to store the variations in \eqref{eq_bfgs_vars_local_2}, since the old variables $ \bbz_i^t$ and $\nabla f_i(\bbz_i^t)$ are available in memory and the updated versions $\bbz_i^{t+1}=\bbx^{t+1}$ and $ \nabla f_i(\bbz_i^{t+1}) = \nabla f_i(\bbx^{t+1}) $ are evaluated at step $t$; see Fig. \ref{fig_sag} for more details. 

Because of the cyclic update scheme, the set of iterates $\{ \bbz_1^t, \bbz_2^t, \hdots, \bbz_n^t\}$ is equal to the set $\{ \bbx^t, \bbx^{t-1}, \hdots, \bbx^{t-n+1}\}$, and, therefore, the set of variables used in the update of IQN is the set of the last $n$ iterates. The update of IQN in \eqref{eq_incr_bfgs_update_222} incorporates the information of all the functions $f_1,\dots,f_n$ to compute the updated variable $\bbx^{t+1}$; however, it uses delayed variables, gradients, and Hessian approximations rather than the the updated variable $\bbx^{t+1}$ for all functions as in classic quasi-Newton methods. The use of delay allows IQN to update the information of a single function at each iteration, thus reducing the computational complexity relative to classic quasi-Newton methods. 

Although the update in \eqref{eq_incr_bfgs_update_222} is helpful in understanding the rationale behind the IQN method, it cannot be implemented at a low computation cost, since it requires computation of the sums $ \sum_{i=1}^n \bbB_i^t $, $\sum_{i=1}^n \bbB_i^t \bbz_i^t $, and $\sum_{i=1}^n \nabla f_i (\bbz_i^t)$ as well as computing the inversion $(\sum_{i=1}^n \bbB_i^t )^{-1}$. In the following section, we introduce an efficient implementation of the IQN method that has the computational complexity of $\mathcal{O}(p^2)$.

\subsection{Efficient implementation of IQN}\label{sec_iqn_implementation} 

To see that the updating scheme in \eqref{eq_incr_bfgs_update_222} requires evaluation of only a single gradient and Hessian approximation matrix per iteration, consider writing the update as
\begin{align}\label{eq_incr_bfgs_update_short}
\bbx^{t+1} =(\tbB^t)^{-1} \left( \bbu^t- \bbg^t \right),
\end{align}
where we define $\tbB^t := \sum_{i=1}^n \bbB_i^t $  as the aggregate Hessian approximation, $\bbu^t := \sum_{i=1}^n \bbB_i^t \bbz_i^t$ as the aggregate Hessian-variable product, and $\bbg^t := \sum_{i=1}^n \nabla f_i (\bbz_i^t)$ as the aggregate gradient. Then, given that at step $t$ only a single index $i_t$ is updated, we can evaluate these variables for step $t+1$ as
\begin{align}
\tbB^{t+1} &= \tbB^t +  \left( \bbB_{i_t}^{t+1} - \bbB_{i_t}^{t} \right), \label{eq_quick_updates_1}\\ 
\bbu^{t+1} &= \bbu^t +\left( \bbB_{i_t}^{t+1} \bbz_{i_t}^{t+1} - \bbB_{i_t}^{t}\bbz_{i_t}^{t} \right),\label{eq_quick_updates_2}
 \\ \bbg^{t+1} &= \bbg^t +  \left(\nabla f_{i_t} (\bbz_{i_t}^{t+1}) - \nabla f_{i_t}(\bbz_{i_t}^{t}) \right). \label{eq_quick_updates_3}
\end{align}
Thus, only $\bbB_{i_t}^{t+1}$ and $\nabla f_{i_t} (\bbz_{i_t}^{t+1})$ are required to be computed at step $t$.

Although the updates in \eqref{eq_quick_updates_1}-\eqref{eq_quick_updates_3} have low computational complexity, the update in \eqref{eq_incr_bfgs_update_short} requires computing $(\tbB^t)^{-1} $ which has a computational complexity of $\ccalO(p^3)$.  This inversion can be avoided by simplifying the update in \eqref{eq_quick_updates_1} as
\begin{align}\label{eq_incr_bfgs_update12}
\tbB^{t+1} =  \tbB^t + \frac{\bby_{i_t}^t \bby^{tT}_{i_t}}{\bby_i^{tT} \bbs^{i_t}_t} - \frac{\bbB_{i_t}^{t} \bbs_{i_t}^t {\bbs^{tT}_{i_t}}\bbB_{i_t}^{t}}{{\bbs^{tT}_{i_t}} \bbB_{i_t}^{t} \bbs_{i_t}^t} .
\end{align}
To derive the expression in \eqref{eq_incr_bfgs_update12} we have substituted the difference $\bbB_{i_t}^{t+1} - \bbB_{i_t}^{t}$ by its rank two expression in \eqref{eq_bfgs_update_local}. Given the matrix $(\tbB^{t})^{-1}$, by applying the Sherman-Morrison formula twice  to the update in \eqref{eq_incr_bfgs_update12} we can compute $(\tbB^{t+1})^{-1}$ as
\begin{align}\label{eq_incr_bfgs_update_inverse}
(\tbB^{t+1})^{-1}&=\bbU^{t} +\frac{\ \!\bbU^{t} (\bbB_{i_t}^{t} \bbs_{i_t}^t)(\bbB_{i_t}^{t} \bbs_{i_t}^t)^T\bbU^{t}}{{\bbs_{i_t}^t}^T\bbB_{i_t}^{t} \bbs_{i_t}^t-(\bbB_{i_t}^{t} \bbs_{i_t}^t)^T\bbU^{t} (\bbB_{i_t}^{t} \bbs_{i_t}^t)},
\end{align}
where the matrix $\bbU^{t}$ is evaluated as 
\begin{align}\label{eq_incr_bfgs_update_inverse_aux}
\bbU^{t} = (\tbB^{t})^{-1} -\frac{(\tbB^{t})^{-1} \bby_{i_t}^t{\bby_{i_t}^{tT}}(\tbB^{t})^{-1}}{\bby_{i_t}^{tT} \bbs_{i_t}^t+{\bby_{i_t}^{tT}}(\tbB^{t})^{-1} \bby_{i_t}^t}.
\end{align}
The computational complexity of the updates in \eqref{eq_incr_bfgs_update_inverse} and \eqref{eq_incr_bfgs_update_inverse_aux} is of the order $\ccalO(p^2)$ rather than the $\ccalO(p^3)$ cost of computing the inverse directly. Therefore, the overall cost of IQN is of the order  $\ccalO(p^2)$ which is substantially lower than $\ccalO(np^2)$ of deterministic quasi-Newton methods.

%
\begin{algorithm}[t]
{\small\begin{algorithmic}[1]
  \REQUIRE $\bbx^{0}$,$\{\nabla f_i (\bbx^0)\}_{i=1}^n$, $\{\bbB_i^0\}_{i=1}^n$
  \STATE Set $\bbz_1^0 =\dots=\bbz_n^0 = \bbx^0$
  \STATE Set ${{(\tbB^0)}}^{-1}=(\sum_{i=1}^n \bbB_i^0)^{-1} $, $\bbu^0 = \sum_{i=1}^n \bbB_i^0 \bbx^0$, $\bbg^0= \sum_{i=1}^n \nabla f_i (\bbx^0)$
  \FOR{$t=0,1,2,\hdots$}
  \STATE Set $i_{t} =(t \mod n) + 1$
  \STATE Compute $\bbx^{t+1} =(\tbB^t)^{-1} \left( \bbu^t- \bbg^t \right)$ [cf. \eqref{eq_incr_bfgs_update_short}]
    \STATE Compute $\bbs_{i_t}^{t+1}$, $\bby_{i_t}^{t+1}$ [cf. \eqref{eq_bfgs_vars_local_2}], and $\bbB_{i_t}^{t+1}$ [cf. \eqref{eq_bfgs_update_local}]
  \STATE Update $\bbu^{t+1}$ [cf. \eqref{eq_quick_updates_2}], $\bbg^{t+1}$ [cf. \eqref{eq_quick_updates_3}], and $(\tbB^{t+1})^{-1}$ [cf. \eqref{eq_incr_bfgs_update_inverse}, \eqref{eq_incr_bfgs_update_inverse_aux}]
   \STATE Update the functions' information tables as in \eqref{eq_local_variable_update_222}, \eqref{eq_local_gradient_update_333}, and \eqref{eq_bfgs_update_local}
  \ENDFOR
\end{algorithmic}}
\caption{Incremental Quasi-Newton (IQN) method}
\label{alg_incr_bfgs}
\end{algorithm}

The complete IQN algorithm is outlined in Algorithm \ref{alg_incr_bfgs}. Beginning with initial variable $\bbx^0$ and gradient and Hessian estimates $\nabla f_i (\bbx^0)$ and $\bbB_i^0$ for all $i$, each variable copy $\bbz^0_i$ is set to $\bbx^0$ in Step 1 and initial values are set for $\bbu^0$, $\bbg^0$ and $(\tbB^0)^{-1}$ in Step 2. For all $t$, in Step 4 the index $i_t$ of the next function to update is selected cyclically. The variable $\bbx^{t+1}$ is computed according to the update in \eqref{eq_incr_bfgs_update_short} in Step 5. In Step 6, the variable $\bbs_{i_t}^{t+1}$ and gradient $\bby_{i_t}^{t+1}$ variations are evaluated as in \eqref{eq_bfgs_vars_local_2} to compute the BFGS matrix $\bbB_{i_t}^{t+1}$ from the update in \eqref{eq_bfgs_update_local}. This information, as well as the updated variable and its gradient, are used in Step 7 to update $\bbu^{t+1}$ and $\bbg^{t+1}$ as in \eqref{eq_quick_updates_2} and \eqref{eq_quick_updates_3}, respectively. The inverse matrix $(\tbB^{t+1})^{-1}$ is also computed by following the expressions in \eqref{eq_incr_bfgs_update_inverse} and \eqref{eq_incr_bfgs_update_inverse_aux}. Finally in Step 8, we update the variable, gradient, and Hessian approximation tables based on the policies in \eqref{eq_local_variable_update_222}, \eqref{eq_local_gradient_update_333}, and \eqref{eq_bfgs_update_local}, respectively.

\section{Convergence Analysis} \label{sec_conv_analysis}
In this section, we study the convergence rate of the proposed IQN method. We first establish its local linear convergence rate, then demonstrate limit properties of the Hessian approximations, and finally show that in a region local to the optimal point the sequence of residuals converges at a superlinear rate. To prove these results we make two main assumptions, both of which are standard in the analysis of quasi-Newton methods. 

\vspace{3mm}
\begin{assumption} \label{ass_gradient_assumption}
There exist positive constants $0 < \mu \leq L$ such that, for all $i$ and $\bbx, \hbx\in \reals^p$, we can write
\begin{equation}\label{eq_gradient_assumption}
\mu \| \bbx - \hbx \|^2 \leq(\nabla f_i (\bbx) - \nabla f_i(\hbx) )^T(\bbx - \hbx) \leq L \| \bbx - \hbx \|^2.
\end{equation}
\end{assumption}
%
\begin{assumption} \label{ass_hessian_assumption}
There exists a positive constant $0 < \tilde{L}$ such that, for all $i$ and $\bbx, \hbx\in \reals^p$, we can write
\begin{equation}\label{eq_hessian_assumption}
 \| \nabla^2 f_i (\bbx) - \nabla^2 f_i(\hbx) \| \leq \tilde{L} \| \bbx - \hbx \|.
\end{equation}
\end{assumption}

The lower bound in \eqref{eq_gradient_assumption} implies that the functions $f_i$ are strongly convex with constant $\mu$, and the upper bound shows that the gradients $\nabla f_i$ are  Lipschitz continuous with parameter $L$. 

The condition in Assumption \ref{ass_hessian_assumption}, states that the Hessians $\nabla^2 f_i$ are Lipschitz continuous with constant $\tilde{L}$. This assumption is commonly made in the analyses of Newton's method \cite{nesterov2013introductory} and quasi-Newton algorithms \cite{Broyden,Powell,Dennis}. According to Lemma~3.1 in \cite{Broyden}, Lipschitz continuity of the Hessians with constant $\tilde{L}$ implies that for $i=1,\dots,n$ and arbitrary vectors $\bbx,\tbx,\hbx\in \reals^p$ we can write 
\begin{equation}\label{result_of_lip_hessian}
\left\| \nabla^2 f_i(\tbx) (\bbx-\hbx) -(\nabla f_i(\bbx)-\nabla f_i(\hbx))\right\|
\leq {\tilde{L}} \|\bbx-\hbx\| \max\left\{\|\bbx-\tbx\|,\|\hbx-\tbx\|\right\}.
\end{equation}
We use the inequality in \eqref{result_of_lip_hessian} in the process of proving the convergence of IQN. 

The goal of BFGS quasi-Newton methods is to approximate the objective function Hessian using the first-order information. Likewise, in the incremental BFGS method, we aim to show that the Hessian approximation matrices for all functions $f_1,\dots,f_n$ are close to the exact Hessian. In the following lemma, we study the difference between the $i$-th optimal Hessian $ \nabla^2 f_i(\bbx^*)$ and its approximation $\bbB_i^t$ over time.  

\vspace{3mm}
\begin{lemma} \label{lemma_M_2}
Consider the proposed IQN method in \eqref{eq_incr_bfgs_update_222}. Further, let $i$ be the index of the updated function at step $t$, i.e., $i=i_t$. 
Define the residual sequence for function $f_i$ as $\sigma_i^t:= \max\{\|\bbz_i^{t+1}-\bbx^*\|,\|\bbz_i^t-\bbx^*\|\}$ and set $\bbM=\nabla^2 f_i(\bbx^*)^{-1/2}$. If Assumptions \ref{ass_gradient_assumption} and \ref{ass_hessian_assumption} hold and the condition $  \sigma_i^t<m/(3\tilde{L})$ is satisfied then
\begin{align}\label{lemma_M_2_claim}
\left\|\bbB_i^{t+1}- \nabla^2 f_i(\bbx^*)\right\|_\bbM \leq \left[ (1- \alpha {\theta_i^t}^2)^{1/2} 
	  +\alpha_3 \sigma_i^t \right]\left\|\bbB_i^t- \nabla^2 f_i(\bbx^*)\right\|_\bbM + \alpha_4 \sigma_i^t,
\end{align}
where $\alpha,\alpha_3 $, and $\alpha_4$ are some positive bounded constants and 
\begin{align}\label{def_theta}
\theta_i^t =\frac{\|\bbM(\bbB_i^t- \nabla^2 f_i(\bbx^*))\bbs_i^t\|}{\|\bbB_i^t- \nabla^2 f_i(\bbx^*)\|_\bbM\|\bbM^{-1}\bbs_i^t\|} \ \ \for \ \bbB_i^t \neq  \nabla^2 f_i(\bbx^*), \quad \theta_i^t=0\ \ \for\  \bbB_i^t=  \nabla^2 f_i(\bbx^*).
\end{align}
\end{lemma}
\begin{proof}
See Appendix \ref{apx_lemma_M_2}.
\end{proof}

\vspace{3mm}

The result in \eqref{lemma_M_2_claim} establishes an upper bound for the weighted norm $\|\bbB_i^{t+1}- \nabla^2 f_i(\bbx^*)\|_\bbM$ with respect to its previous value $\|\bbB_i^{t}- \nabla^2 f_i(\bbx^*)\|_\bbM$ and the sequence $\sigma_i^t:= \max\{\|\bbz_i^{t+1}-\bbx^*\|,\|\bbz_i^t-\bbx^*\|\}$, when the variables are in a neighborhood of the optimal solution such that $  \sigma_i^t<m/(3\tilde{L})$. Indeed, the result in  \eqref{lemma_M_2_claim} holds only for the index $i=i_t$ and for the rest of indices we have $\|\bbB_i^{t+1}- \nabla^2 f_i(\bbx^*)\|_\bbM=\|\bbB_i^{t}- \nabla^2 f_i(\bbx^*)\|_\bbM$ simply by definition of the cyclic update. Note that if the residual sequence $\sigma_i^t$ associated with $f_i$ approaches zero, we can simplify \eqref{lemma_M_2_claim} as  
\begin{align}\label{bfgs_dec_ineq}
\|\bbB_i^{t+1}- \nabla^2 f_i(\bbx^*)\|_\bbM \lesssim (1- \alpha {\theta_i^t}^2)^{1/2} 
	 \|\bbB_i^t- \nabla^2 f_i(\bbx^*)\|_\bbM.
	 \end{align} 
The equation in \eqref{bfgs_dec_ineq} implies that if ${\theta_i^t}$ is always strictly larger than zero, the sequence $\|\bbB_i^{t+1}- \nabla^2 f_i(\bbx^*)\|_\bbM $ approaches zero. If not, then the sequence ${\theta_i^t}$ converges to zero which implies the  Dennis-Mor\'e condition from \eqref{eq_dennis_more}, i.e.
\begin{equation}\label{bfgs_dec_ineq_2}
\lim_{t \to \infty} \frac{\|(\bbB_i^t- \nabla^2 f_i(\bbx^*))\bbs_i^t\|}{\|\bbs_i^t\|}=0.
\end{equation}
Therefore, under both conditions the result in \eqref{bfgs_dec_ineq_2} holds. This is true since the limit $\lim_{t\to \infty }\|\bbB_i^{t+1}- \nabla^2 f_i(\bbx^*)\|_\bbM =0 $ yields the result in \eqref{bfgs_dec_ineq_2}. 

Based on this intuition, we proceed to show that the sequence $\sigma_i^t$ converges to zero for all $i=1,\dots,n$. To do so, we show that the sequence $\|\bbz_i^t-\bbx^*\|$ is linearly convergent for all $i=1,\dots,n$. To achieve this goal we first prove an upper bound for the error $\|\bbx^{t+1}-\bbx^*\|$ of IQN in the following lemma.

\begin{lemma} \label{lemma_cyclic_property}
Consider the proposed IQN method in \eqref{eq_incr_bfgs_update_222}. If the conditions in Assumptions \ref{ass_gradient_assumption} and \ref{ass_hessian_assumption} hold, then the sequence of iterates generated by IQN satisfies 
\begin{align}\label{lemma_claim_imp20}
&\|\bbx^{t+1}-\bbx^*\|
\leq
 \frac{\tilde{L}\Gamma^{t}}{n} \sum_{i=1}^n  \left\|\bbz_i^t-\bbx^* \right\|^2
 +
 \frac{\Gamma^{t}}{n} \sum_{i=1}^n  \left\| \left(\bbB_i^t-\nabla^2 f_{i}(\bbx^*)\right) \left(\bbz_i^t-\bbx^*\right) \right\|, 
\end{align}
where $\Gamma^t:=\|( (1/n)\sum_{i=1}^n\bbB_i^t )^{-1}\|$.
\end{lemma}

\vspace{3mm}

\begin{proof}
See Appendix \ref{apx_lemma_cyclic_property}.
\end{proof}

\vspace{3mm}

Lemma \ref{lemma_cyclic_property} shows that the residual $\|\bbx^{t+1}-\bbx^*\|$ is bounded above by a sum of quadratic and linear terms of the last $n$ residuals. This can eventually lead to a superlinear convergence rate by establishing the linear term converges to zero at a fast rate, leaving us with an upper bound of quadratic terms only. First, however, we establish a local linear convergence rate in the proceeding theorem to show that the sequence $\sigma_i^t$ converges to zero.

\vspace{3mm}
\begin{lemma}\label{theorem_local_lin_conv}
Consider the proposed IQN method in \eqref{eq_incr_bfgs_update_222}. If Assumptions \ref{ass_gradient_assumption} and \ref{ass_hessian_assumption} hold, then, for any $r\in(0,1)$ there are positive constants $\eps(r)$ and $\delta(r)$ such that if $\|\bbx^0-\bbx^*\|< \eps(r)$ and $\|\bbB_i^0-\nabla^2 f_i(\bbx^*)\|_\bbM< \delta(r)$ for $\bbM=\nabla^2 f_i(\bbx^*)^{-1/2}$ and $i=1,2,\dots,n$, the sequence of iterates generated by IQN satisfies 
\begin{equation} \label{loc_lin_convg55}
\|\bbx^{t}-\bbx^*\|\leq r^{[\frac{t-1}{n}]+1} \|\bbx^0-\bbx^*\|.
\end{equation}
Moreover, the sequences of norms $\{\|\bbB_i^t\|\}$ and $\{\|({\bbB_i^t})^{-1}\|\}$ are uniformly bounded. 
\end{lemma}

\vspace{3mm}
\begin{proof}
See Appendix \ref{apx_theorem_local_lin_conv}.
\end{proof}

\vspace{3mm}

The result in Lemma \ref{theorem_local_lin_conv} shows that the sequence of iterates generated by IQN has a local linear convergence rate after each pass over all functions. Consequently, we obtain that the $i$-th residual sequence $\sigma_i^t$ is linearly convergent for all $i$. Note that Lemma \ref{theorem_local_lin_conv} can be considered as an extension of Theorem 3.2 in \cite{Broyden} for incremental settings.
Following the arguments in \eqref{bfgs_dec_ineq} and \eqref{bfgs_dec_ineq_2}, we use the summability of the sequence $\sigma_i^t$ along with the result in Lemma \ref{lemma_M_2} to prove Dennis-Mor\'e condition for all functions $f_i$. 

\vspace{3mm}
\begin{proposition}\label{prop_DM_condition}
Consider the proposed IQN method in \eqref{eq_incr_bfgs_update_222}. Assume that the hypotheses in Lemmata \ref{lemma_M_2} and \ref{theorem_local_lin_conv} are satisfied. Then, for all $i=1,\dots,n$ it holds,
\begin{equation}\label{instan_dennis_more}
\lim_{t \to \infty} \frac{\|(\bbB_i^t- \nabla^2 f_i(\bbx^*))\bbs_i^t\|}{\|\bbs_i^t\|}=0.
\end{equation}
\end{proposition}

\vspace{3mm}
\begin{proof}
See Appendix \ref{apx_prop_DM_condition}.
\end{proof}
\vspace{3mm}

The statement in Proposition \ref{prop_DM_condition} indicates that for each function $f_i$ the Dennis-Mor\'e condition holds. In the tradition quasi-Newton methods the Dennis-Mor\'e condition is sufficient to show that the method is superlinearly convergent. However, the same argument does not hold for the proposed IQN method, since we can't recover the Dennis-Mor\'e condition for the global objective function $f$ from the result in Proposition \ref{prop_DM_condition}. In other words, the result in \eqref{instan_dennis_more} does not imply the limit in \eqref{eq_dennis_more} required in the superlinear  convergence analysis of quasi-Newton methods. Therefore, here we pursue a different approach and seek to prove that the linear terms $(\bbB_i^t-\nabla^2 f_{i}(\bbx^*)) (\bbz_i^t-\bbx^*) $ in \eqref{lemma_claim_imp20} converge to zero at a superlinear rate, i.e., for all $i$ we can write $lim_{t\to\infty} \|(\bbB_i^t-\nabla^2 f_{i}(\bbx^*)) (\bbz_i^t-\bbx^*) \|/\|\bbz_i^t-\bbx^*\|=0$. If we establish this result, it follows from the result in Lemma \ref{lemma_cyclic_property} that the sequence of residuals $\|\bbx^t-\bbx^*\|$ converges to zero superlinearly. 

We continue the analysis of the proposed IQN method by establishing a generalized limit property that follows from the Dennis-Mor\'e criterion in \eqref{eq_dennis_more}. In the following lemma, we leverage  the local linear convergence of the iterates $\bbx^t$ to show that that the vector $\bbz_i^t-\bbx^*$ lies in the null space of $\bbB_i^t - \nabla^2 f_i(\bbx^*)$ as $t$ approaches infinity.

\vspace{3mm}

\begin{lemma} \label{lemma_incr_bfgs_lim_final}
Consider the proposed IQN method in \eqref{eq_incr_bfgs_update_222}. Assume that the hypotheses in Lemmata \ref{lemma_M_2} and \ref{theorem_local_lin_conv} are satisfied.  As $t$ approaches infinity, the following holds for all $i$,
\begin{equation}\label{eq_lemma_incr_bfgs_lim_final}
\lim_{t \rightarrow \infty} \frac{ \| (\bbB_i^t - \nabla^2 f_i(\bbx^*)) (\bbz_i^t - \bbx^*) \|}{\| \bbz_i^t -\bbx^*\|} = 0.
\end{equation}
\end{lemma}
\vspace{3mm}
\begin{proof}
See Appendix \ref{apx_lemma_incr_bfgs_lim_final}.
\end{proof}
\vspace{3mm}

The result in Lemma \ref{lemma_incr_bfgs_lim_final} can thus be used in conjunction with Lemma \ref{lemma_cyclic_property}  to show that the residual $\|\bbx^{t+1} - \bbx^*\|$ is bounded by a sum of quadratic terms of previous residuals and a term that converges to zero superlinearly. This result leads us to the following result, namely the local superlinear convergence of the sequence of residuals with respect to the average sequence, stated in the following theorem.

\vspace{3mm}
\begin{theorem}\label{sup_linear_thm}
Consider the proposed IQN method in \eqref{eq_incr_bfgs_update_222}. Suppose that the conditions in the hypotheses of Lemmata \ref{lemma_M_2} and \ref{theorem_local_lin_conv} are valid. Then, the sequence of residuals $ \| \bbx^t - \bbx^*\|$ satisfies 
\begin{equation}\label{final_claim}
\lim_{t \rightarrow \infty} \frac{ \| \bbx^{t} - \bbx^*\|}{\frac{1}{n} (\| \bbx^{t-1} - \bbx^*\|+\dots+ \| \bbx^{t-n} - \bbx^*\|)} = 0.
\end{equation}
\end{theorem}
\vspace{3mm}

\begin{proof}
See Appendix \ref{apx_sup_linear_thm}.
\end{proof}
\vspace{3mm}
The result in \eqref{final_claim} shows a mean-superlinear convergence rate for the sequence of iterates generated by IQN. To be more precise, it shows that the ratio that captures the error at step $t$ divided by the average of last $n$ errors converges to zero. This is not equivalent to the classic Q-superlinear convergence for full-batch quasi-Newton methods, i.e., $\lim_{t\to\infty} \|\bbx^{t+1}-\bbx^*\|/\|\bbx^{t}-\bbx^*\|=0$. Although Q-superlinear convergence of the residuals $\|\bbx^{t}-\bbx^*\|$ is not provable, we can show that there exists a subsequence of the sequence $\|\bbx^{t}-\bbx^*\|$ that converges to zero superlinearly. In addition, there exists a superlinearly convergent sequence that is an upper bound for the original sequence of errors $\|\bbx^{t}-\bbx^*\|$. We formalize these results in the following theorem.

\vspace{3mm}
\begin{theorem}\label{sup_linear_thm_2}
Consider the proposed IQN method in \eqref{eq_incr_bfgs_update_222}. Suppose that the conditions in the hypotheses of Lemmata \ref{lemma_M_2} and \ref{theorem_local_lin_conv} are valid. Then, there exists a subsequence of $\|\bbx^{t}-\bbx^*\|$ that converges to zero superlinearly. Moreover, there exists a sequence $\zeta^t$ such that $\|\bbx^t-\bbx^*\|\leq \zeta^t$ for all $t\geq0$, and the sequence $\zeta^t$ converges to zero at a superlinear rate, i.e., 
\begin{equation}\label{sup_linear_thm_2_claim}
\lim_{t\to\infty} \frac{\zeta^{t+1}}{\zeta^t}=0.
\end{equation}
\end{theorem}
\begin{proof}
See Appendix \ref{apx_sup_linear_thm_2}.
\end{proof}
\vspace{3mm}

The first result in Theorem \ref{sup_linear_thm_2} states that although the whole sequence $\|\bbx^{t}-\bbx^*\|$ is not necessarily superlinearly convergent, there exists a subsequence of the sequence $\|\bbx^{t}-\bbx^*\|$ that converges at a superlinear rate. The second claim in Theorem \ref{sup_linear_thm_2} establishes R-superlinear convergence rate of the whole sequence $\|\bbx^{t}-\bbx^*\|$. In other words, it guarantees that $\|\bbx^{t}-\bbx^*\|$ is upper bounded by a superlinearly convergent sequence.

\section{Related Works} \label{sec_related_works}

Various methods have been studied in the literature to improve the performance of traditional full-batch optimization algorithms. The most famous method for reducing the computational complexity of gradient descent (GD) is stochastic gradient descent (SGD), which uses the gradient of a single randomly chosen function to approximate the full-gradient \cite{Bottou}.  Incremental gradient descent method (IGD) is similar to SGD except the function is chosen in a cyclic routine \cite{blatt2007convergent}. Both SGD and IGD suffer from slow sublinear convergence rate because of the noise of gradient approximation. The incremental aggregated methods, which use memory to aggregate the gradients of all $n$ functions, are successful in reducing the noise of gradient approximation to achieve linear convergence rate \cite{roux2012stochastic,Schmidt2016,defazio2014saga,johnson2013accelerating}. The work in \cite{roux2012stochastic} suggests a random selection of functions which leads to stochastic average gradient method (SAG), while the works in \cite{blatt2007convergent,gurbuzbalaban2015convergence,mokhtari2016surpassing} use a cyclic scheme.


Moving beyond first order information, there have been stochastic quasi-Newton methods to approximate Hessian information \cite{SchraudolphYG07,mokhtari2014res,mokhtari2015global,MoritzNJ16,gower2016stochastic}. All of these stochastic quasi-Newton methods reduce computational cost of quasi-Newton methods by updating only a randomly chosen single or small subset of gradients at each iteration. However, they are not able to recover the superlinear convergence rate of quasi-Newton methods \cite{Broyden,Powell,Dennis}. The incremental Newton method (NIM) in \cite{rodomanov2016superlinearly} is the only incremental method shown to have a superlinear convergence rate; however, the Hessian function is not always available or computationally feasible. Moreover, the implementation of NIM requires computation of the incremental aggregated Hessian inverse which has the computational complexity of the order $\mathcal{O}(p^3)$.

\section{Numerical Results} \label{sec_num_results}

We proceed by simulating the performance of IQN on a variety of machine learning problems on both artificial and real datasets. We compare the performance of IQN against a collection of well known first order stochastic and incremental algorithms---namely SAG, SAGA, and IAG. To begin, we look at a simple quadratic program, also equivalent to the solution of linear least squares estimation problem. Consider the objective function to be minimized,
\begin{align}\label{eq_quadratic_problem}
\bbx^* = \argmin_{\bbx\in\reals^p} f(\bbx) := \argmin_{\bbx\in\reals^p} \frac{1}{n}\sum_{i=1}^n \frac{1}{2} \bbx^T \mathbf{A}_i \bbx + \mathbf{b}_i^T \bbx.
\end{align}
We generate $\mathbf{A}_i \in \reals^{p \times p}$ as a random positive definite matrix and $\mathbf{b}_i \in \reals^p$ as a random vector for all $i$. In particular we set the matrices $\mathbf{A}_i := \text{diag}\{ \mathbf{a}_i \}$ and generate random vectors $\bba_i$ with the first $p/2$ elements chosen from $[1, 10^{\xi/2}]$ and last $p/2$ elements chosen from $[10^{-\xi/2}, 1]$.  The parameter $\xi$ is used to manually set the condition number for the quadratic program in \eqref{eq_quadratic_problem}, ranging from $\xi = 1$ (i.e. small condition number $10^2$) and $\xi=2$ (i.e. large condition number $10^4$). The vectors $\mathbf{b}_i$ are chosen uniformly and randomly from the box $[0,10^3]^p$. The variable dimension is set to be $p=10$ and number of functions $n=1000$. Given that we focus on local convergence, we use a constant step size of $\eta = 1$ for the proposed IQN method while choosing the largest step size allowable by the other methods to converge.

\begin{figure}[t]
\centering
\includegraphics[width=0.45\linewidth,height=0.3\linewidth]{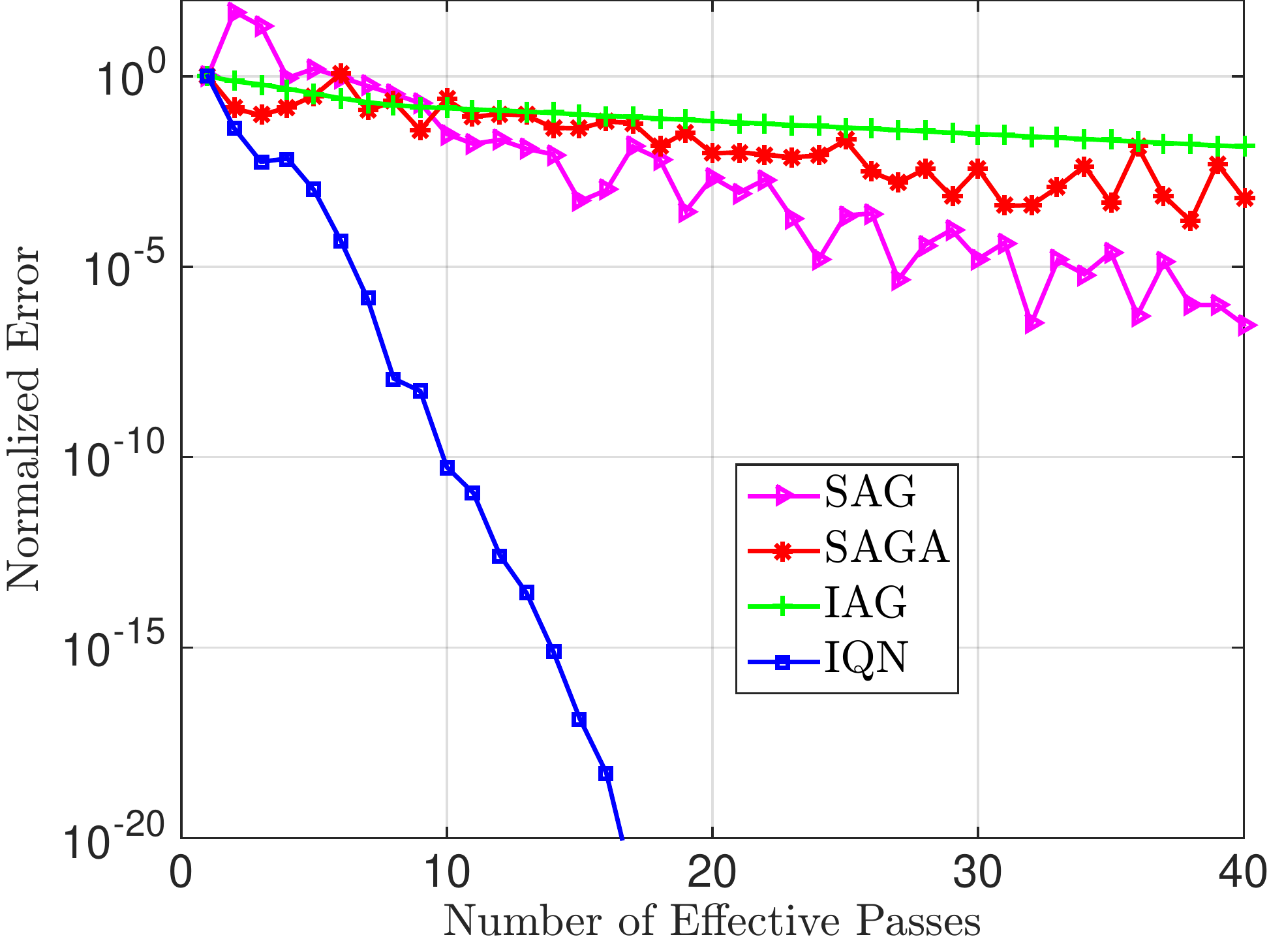}\qquad
\includegraphics[width=0.45\linewidth,height=0.3\linewidth]{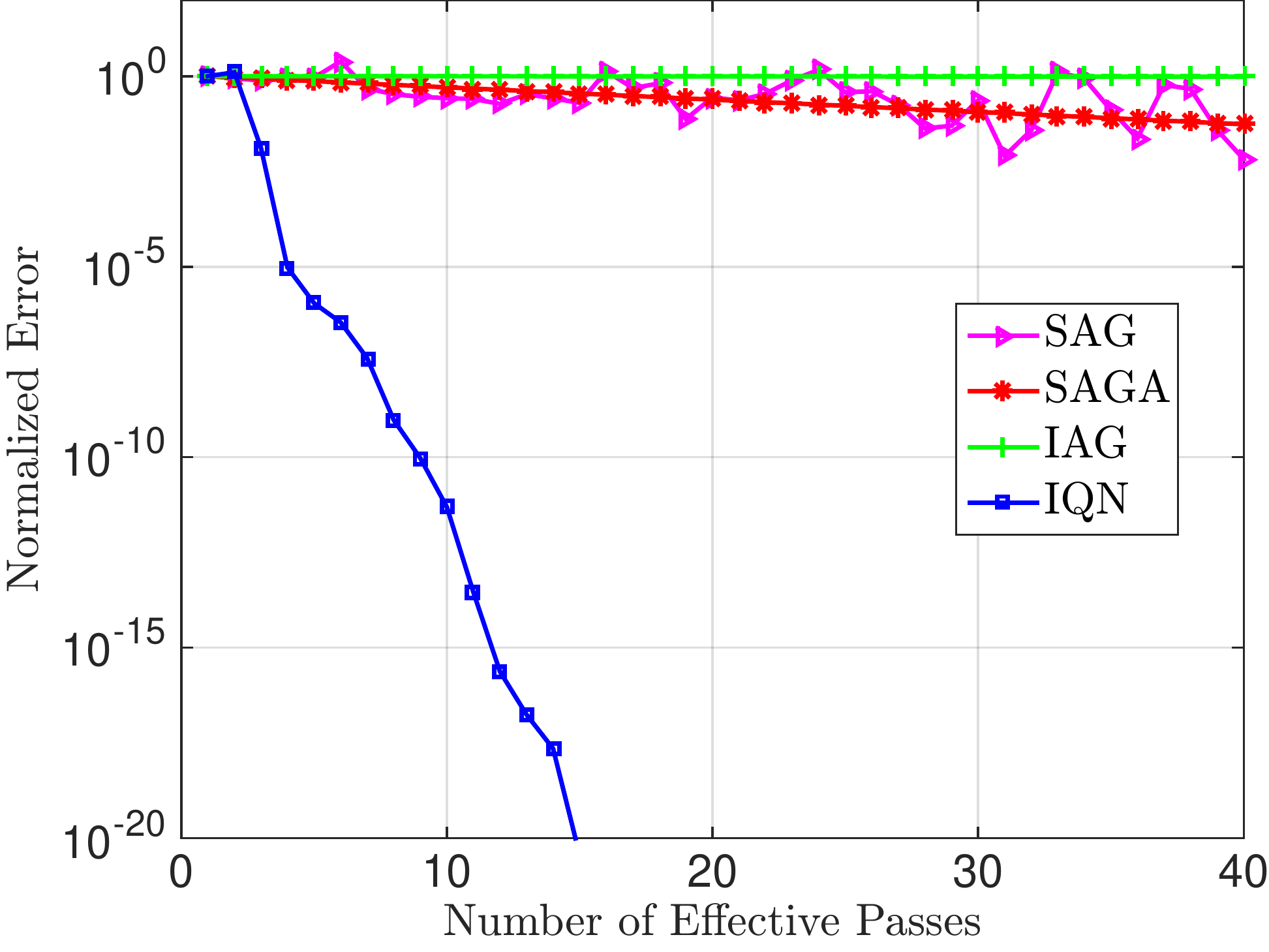} 
\caption{Convergence results of proposed IQN method in comparison to SAG, SAGA, and IAG. In the left image, we present a sample convergence path of the normalized error on the quadratic program with a small condition number. In the right image, we show the convergence path for the quadratic program with a large condition number. In all cases, IQN provides significant improvement over first order methods, with the difference increasing for larger condition number.}\label{fig_quad_results}
\end{figure}

In Figure \ref{fig_quad_results} we present a simulation of the convergence path of the normalized error $\|\bbx^t-\bbx^*\|/\|\bbx^0-\bbx^*\|$ for the quadratic program. In the the left image, we show a sample simulation path for all methods on the quadratic problem with a small condition number. Step sizes of $\eta=5 \times 10^{-5}$, $\eta=10^{-4}$ and $\eta=10^{-6}$ were used for SAG, SAGA, and IAG, respectively. These step sizes are tuned to compare the best performance of these methods with IQN. The proposed method reaches a error of $10^{-10}$ after 10 passes through the data. Alternatively, SAGA achieves the same error of $10^{-5}$ after 30 passes, while SAG and IAG do not reach $10^{-5}$ after 40 passes. 

In the right image of Figure \ref{fig_quad_results}, we repeat the same simulation but with larger condition number. In this case, SAG uses stepsize $\eta=2\times10^{-4}$ while others remain the same. Observe that while the performance of IQN does not degrade with larger condition number, the first order methods all suffer large degradation. SAG, SAGA, and IAG reach after 40 passes a normalized error of $6.5 \times 10^{-3}$, $5.5 \times 10^{-2}$, and $9.6 \times 10^{-1}$, respectively. It can be seen that IQN significantly outperforms the first order method for both condition number sizes, with the outperformance increasing for larger condition number. This is an expected result, as first order methods often do not perform well for ill conditioned problems. 

\begin{figure}[t]
\centering
\includegraphics[width=0.45\linewidth,height=0.3\linewidth]{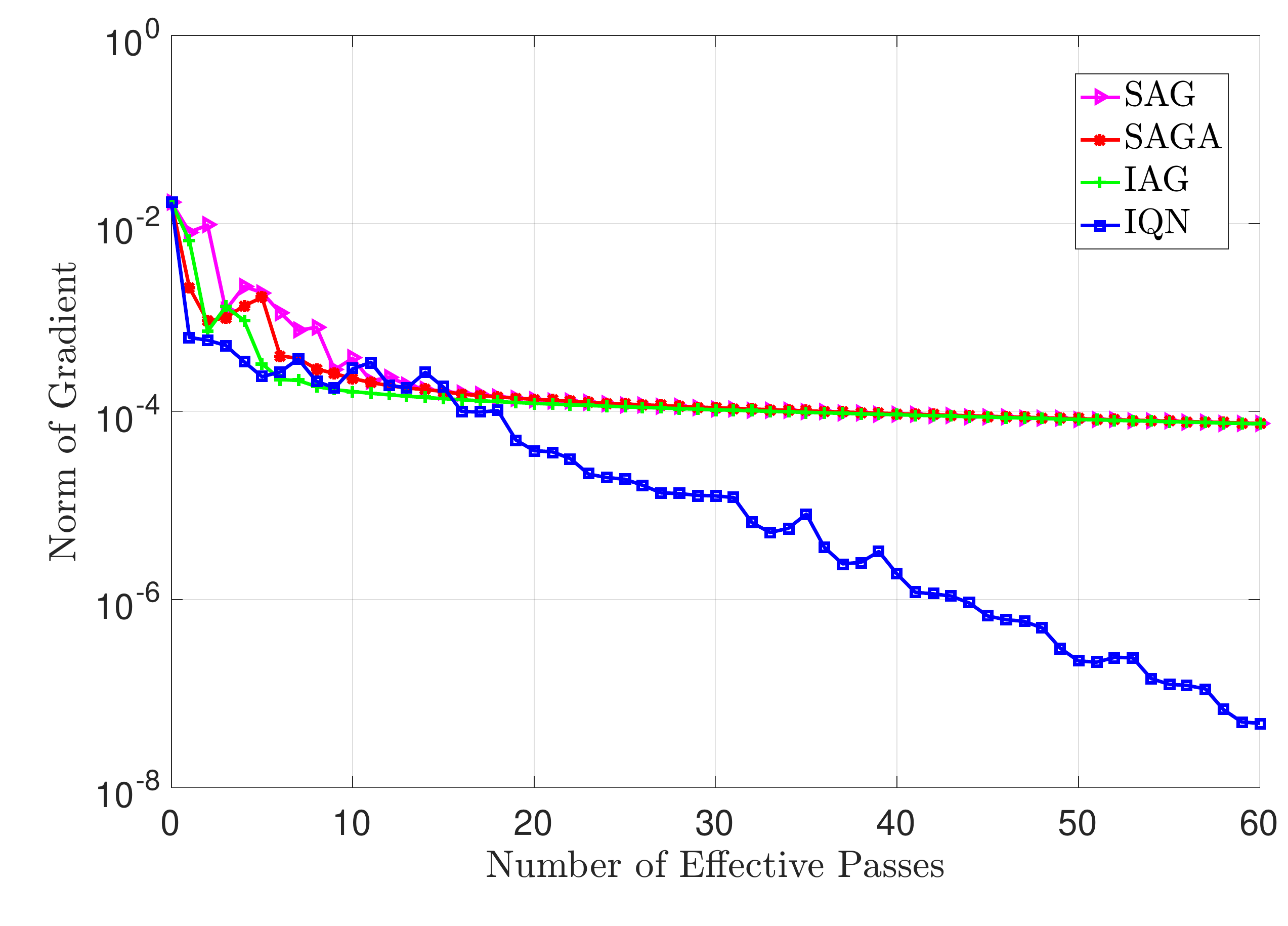}
\caption{Convergence results for a sample convergence path for the logistic regression problem on classifying handwritten digits. IQN substantially outperforms the first order methods.}\label{fig_log_results}
\end{figure}

\subsection{Logistic regression}

We proceed to numerically evaluate the performance of IQN relative to existing methods on the classification of handwritten digits in the MNIST database \cite{lecun2010mnist}. In particular, we solve the binary logistic regression problem.  A logistic regression takes as inputs $n$ training feature vectors $\bbu_i \in \reals^p$ with associated labels $v_i \in \{-1,1\}$ and outputs a linear classifier $\bbx$ to predict the label of unknown feature vectors. For the digit classification problem, each feature vector $\bbu_i$ represents a vectorized image and label $v_i$ its label as one of two digits. We evaluate for any training sample $i$ the probability of a label $v_i=1$ given image $\bbu_i$ as $P(v=1 | \bbu) = 1/(1+\exp(-\bbu^T\bbx))$. The classifier $\bbx$ is chosen to be the vector which maximizes the log likelihood across all $n$ samples. Given $n$ images $\bbu_i$ with associated labels $v_i$, the optimization problem for logistic regression is written as
\begin{align}
\bbx^* = \argmin_{\bbx\in\reals^p} f(\bbx) := \argmin_{\bbx\in\reals^p}\frac{\lambda}{2}\| \bbx\|^2 + \frac{1}{n}\sum_{i=1}^{n} \log[ 1 + \exp(-v_{i}\bbu_{i}^T\bbx)],
\label{eq_logistic_problem}
\end{align}
where the first term is a regularization term parametrized by $\lambda \geq 0$. 

For our simulations we select from the MNIST dataset $n=1000$ images with dimension $p=784$ labelled as one of the digits ``0" or ``8' and fix the regularization parameter as $\lambda= 1/n$ and stepsize $\eta=0.01$ for all first order methods. In Figure \ref{fig_log_results} we present the convergence path of IQN relative to existing methods in terms of the norm of the gradient. As in the case of the quadratic program, the IQN performs all gradient-based methods. IQN reaches a gradient magnitude of $4.8 \times 10^{-8}$ after 60 passes through the data while the SAGA reaches only a magnitude of $7.4 \times 10^{-5}$ (all other methods perform even worse). Further note that while the first order methods begin to level out after 60 passes, the IQN method continues to descend. These results demonstrate the effectiveness of IQN on a practical machine learning problem with real world data.


%
%

\appendix

\section{Proof of Lemma \ref{lemma_M_2}}\label{apx_lemma_M_2} 
To prove the claim in Lemma \ref{lemma_M_2}, we first prove the the following lemma which is based on the result in \cite[Lemma 5.2]{Broyden}.

\begin{lemma} \label{lemma_M_1}
Consider the proposed IQN method in \eqref{eq_incr_bfgs_update_222}.
Let $\bbM$ be a nonsingular symmetric matrix such that 
\begin{equation}\label{ness_cond}
\|\bbM\bby_i^t-\bbM^{-1}\bbs_i^t\| \leq \beta \|\bbM^{-1}\bbs_i^t\| ,
\end{equation}
for some $\beta \in[0,1/3]$ and vectors $\bbs_i^t$ and $\bby_i^t$ in $\reals^p$ with $\bbs_i^t\neq\bb0$. Consider $i$ as the index of the updated function at step $t$, i.e., $i=i_t$, and let $\bbB_i^t$ be symmetric and computed according to the update in \eqref{eq_bfgs_update_local}. Then, there exist positive constants $\alpha$, $\alpha_1$, and $\alpha_2$ such that, for any symmetric $\bbA\in \reals^{p\times p}$ we have, 
\begin{align}\label{B_cons_rela}
\|\bbB_i^{t+n}-\bbA\|_\bbM \leq \left[ (1- \alpha \theta^2)^{1/2} 
	+\alpha_1 \frac{\|\bbM\bby_i^t-\bbM^{-1}\bbs_i^t\| }{\|\bbM^{-1}\bbs_i^t\|} \right] \|\bbB_i^t-\bbA\|_\bbM + \alpha_2 \frac{\|\bby_i^t-\bbA\bbs_i^t\|}{\|\bbM^{-1}\bbs_i^t\|},
\end{align}
where $\alpha= (1-2\beta)/(1-\beta^2)\in [3/8,1]$, $\alpha_1=2.5(1-\beta)^{-1}$, $\alpha_2=2(1+2\sqrt{p})\|\bbM\|_\bbF$, and
\begin{align}
\theta =\frac{\|\bbM(\bbB_i^t-\bbA)\bbs_i^t\|}{\|\bbB_i^t-\bbA\|_\bbM\|\bbM^{-1}\bbs_i^t\|} \quad \for \ \bbB_i^t\neq \bbA, \qquad  \theta=0\quad  \for\   \bbB_i^t = \bbA.
\end{align}
\end{lemma}

\begin{proof}
First note that the Hessian approximation $\bbB_i^{t+n}$ is equal to $\bbB_i^{t+1}$ if the function $f_i$ is updated at step $t$. Considering this observation and the result of Lemma 5.2. in \cite{Broyden} the claim in \eqref{B_cons_rela} follows.
\end{proof}

The result in Lemma \ref{lemma_M_1} provides an upper bound for the difference between the Hessian approximation matrix $\bbB_i^{t+n}$ and any positive definite matrix $\bbA$ with respect to the difference between the previous Hessian approximation $\bbB_i^{t}$ and the matrix $\bbA$. The interesting choice for the arbitrary matrix $\bbA$ is the Hessian of the $i$-th function at the optimal argument, i.e., $\bbA=\nabla^2 f_i(\bbx^*)$, which allows us to capture the difference between the sequence of Hessian approximation matrices for function $f_i$ and the Hessian $\nabla^2 f_i(\bbx^*)$ at the optimal argument. 
We proceed to use the result in Lemma \ref{lemma_M_1} for $\bbM=\nabla^2f_i(\bbx^*)^{-1/2}$ and $\bbA=\nabla^2f_i(\bbx^*)$ to prove the claim in \eqref{lemma_M_2_claim}. To do so, we first need to show that the condition in \eqref{ness_cond} is satisfied. Note that according to the condition in Assumptions \ref{ass_gradient_assumption} and \ref{ass_hessian_assumption} we can write
\begin{align}\label{proof10}
 \frac{\| \bby_i^t-\nabla^2f_i(\bbx^*)\bbs_i^t\|}{\|\nabla^2f_i(\bbx^*)^{1/2}\bbs_i^t\|}
\leq \frac{ \tilde{L}\|\bbs_i^t\|\max\{\|\bbz_i^t-\bbx^*\|,\|\bbz_i^{t+1}-\bbx^*\|\} }{\sqrt{m}\|\bbs_i^t\|} = \frac{\tilde{L}}{\sqrt{m}} \sigma_i^t
\end{align}
This observation implies that the left hand side of the condition in \eqref{ness_cond} for $\bbM=\nabla^2f_i(\bbx^*)^{-1/2}$ is bounded above by 
\begin{align}\label{proof20}
\frac{\|\bbM\bby_i^t-\bbM^{-1}\bbs_i^t\|}{\|\bbM^{-1}\bbs_i^t\|}
\leq
\frac{\| \nabla^2f_i(\bbx^*)^{-1/2}\| \|\bby_i^t-\nabla^2f_i(\bbx^*)\bbs_i^t\|}{\|\nabla^2f_i(\bbx^*)^{1/2}\bbs_i^t\|}
 \leq \frac{\tilde{L}}{m} \sigma_i^t
\end{align}
Thus, the condition in \eqref{ness_cond} is satisfied since $\tilde{L} \sigma_i^t/m<1/3.$ Replacing the upper bounds in \eqref{proof10} and \eqref{proof20} into the expression in \eqref{B_cons_rela} implies the claim in \eqref{lemma_M_2_claim} with
\begin{align}\label{proof300}
\beta= \frac{\tilde{L}}{m} \sigma_i^t, \ 
 \alpha=\frac{1-2\beta}{1-\beta^2}, \ 
\alpha_3=\frac{5\tilde{L}}{2m(1-\beta)}, \ 
\alpha_4=    \frac{{2(1+2\sqrt{p}) \tilde{L}}}{\sqrt{m}}\|\nabla^2f_i(\bbx^*)^{-\frac{1}{2}}\|_\bbF,
\end{align}
and the proof is complete.

\section{Proof of Lemma \ref{lemma_cyclic_property}}\label{apx_lemma_cyclic_property}
 Start by subtracting $\bbx^*$ from both sides of \eqref{eq_incr_bfgs_update_222} to obtain
\begin{align}\label{lemma_claim_imp_21220}
\bbx^{t+1}-\bbx^*
&=\left( \frac{1}{n}\sum_{i=1}^n\bbB_i^t \right)^{-1}
			 \left( \frac{1}{n} \sum_{i=1}^n\bbB_i^t \bbz_i^t - \frac{1}{n} \sum_{i=1}^n \nabla f_i (\bbz_i^t) - \frac{1}{n} \sum_{i=1}^n\bbB_i^t\bbx^* \right).
\end{align}
As the gradient of $f$ at the optimal point is the vector zero, i.e., $(1/n) \sum_{i=1}^n \nabla f_i (\bbx^*)=\mathbf{0}$, we can subtract $(1/n) \sum_{i=1}^n \nabla f_i (\bbx^*)$ from the right hand side of \eqref{lemma_claim_imp_21220} and rearrange terms to obtain
\begin{equation}\label{eq_3a}
\bbx^{t+1}-\bbx^*=\left( \frac{1}{n}\sum_{i=1}^n\bbB_i^t \right)^{-1}\!\!
 \left( \frac{1}{n} \sum_{i=1}^n\bbB_i^t \left(\bbz_i^t-\bbx^*\right) - \frac{1}{n} \sum_{i=1}^n \left(\nabla f_i (\bbz_i^t)-\nabla f_i (\bbx^*)\right)  \right)\!.
\end{equation}
The expression in \eqref{eq_3a} relates the residual at time $t+1$ to the previous $n$ residuals and the Hessian approximations $\bbB_i^t$. To analyze this further, we can replace the Hessian approximations $\bbB_i^t$ with the actual Hessians $\nabla^2 f_{i}(\bbx^*)$ and the approximation difference $\nabla^2 f_{i}(\bbx^*) - \bbB_i^t$. To do so, we add and subtract $ ({1}/{n}) \sum_{i=1}^n\nabla^2 f_{i}(\bbx^*) \left(\bbz_i^t-\bbx^*\right)$ to the right hand side of \eqref{eq_3a} and rearrange terms to obtain
\begin{align}
\bbx^{t+1}-\bbx^*
&
=\left( \frac{1}{n}\sum_{i=1}^n\bbB_i^t \right)^{-1}
 \left( \frac{1}{n} \sum_{i=1}^n \left[\nabla^2 f_{i}(\bbx^*) \left(\bbz_i^t-\bbx^*\right)-  \left(\nabla f_i (\bbz_i^t)-\nabla f_i (\bbx^*)\right) \right] \right) \nonumber\\
 &\qquad 
 +\left( \frac{1}{n}\sum_{i=1}^n\bbB_i^t \right)^{-1} \left( \frac{1}{n} \sum_{i=1}^n  \left[\bbB_i^t-\nabla^2 f_{i}(\bbx^*)\right] \left(\bbz_i^t-\bbx^*\right) \right).
\end{align}
We proceed to take the norms of both sides and use the triangle inequality to obtain an upper bound on the norm of the residual $\|\bbx^{t+1}-\bbx^*\|$, 
\begin{align}\label{eqr}
\|\bbx^{t+1}-\bbx^*\|
  &
\leq\left\|\left( \frac{1}{n}\sum_{i=1}^n\bbB_i^t \right)^{-1}\right\|
 \frac{1}{n} \sum_{i=1}^n  \left\|\nabla^2 f_{i}(\bbx^*) \left(\bbz_i^t-\bbx^*\right)-  \left(\nabla f_i (\bbz_i^t)-\nabla f_i (\bbx^*)\right)  \right\| \nonumber\\
 &\qquad 
 +\left\|\left( \frac{1}{n}\sum_{i=1}^n\bbB_i^t \right)^{-1}\right\| 
 \frac{1}{n} \sum_{i=1}^n  \left\| \left[\bbB_i^t-\nabla^2 f_{i}(\bbx^*)\right] \left(\bbz_i^t-\bbx^*\right) \right\|.
\end{align}
To obtain the quadratic term in \eqref{lemma_claim_imp20} from the first term in \eqref{eqr}, we use the Lipschitz continuity of the Hessians $\nabla^2 f_{i}$ which leads to the inequality 
\begin{align}\label{eqr_new}
\|\nabla^2 f_{i}(\bbx^*) \left(\bbz_i^t-\bbx^*\right)-  \left(\nabla f_i (\bbz_i^t)-\nabla f_i (\bbx^*)\right) \|
\leq {\tilde{L}}\left\|\bbz_i^t-\bbx^*\right\|^2.
\end{align}
Replacing the expression $\|\nabla^2 f_{i}(\bbx^*) \left(\bbz_i^t-\bbx^*\right)-  \left(\nabla f_i (\bbz_i^t)-\nabla f_i (\bbx^*)\right)\|$ in \eqref{eqr} by the upper bound in \eqref{eqr_new}, the claim in \eqref{lemma_claim_imp20} follows.

\section{Proof of Lemma \ref{theorem_local_lin_conv}}\label{apx_theorem_local_lin_conv}
In this proof we use some steps in the proof of \cite[Theorem 3.2]{Broyden}. 
To start we use the fact that in a finite-dimensional vector space there always  exists a constant $\eta>0$ such that $\|\bbA\|\leq \eta \|\bbA\|_\bbM$. Consider $\gamma=1/m$ is an upper bound for the norm $\|\nabla^2 f(\bbx^*)^{-1}\|$. Assume that $\eps(r)=\eps$ and $\delta(r)=\delta$ are chosen such that 
\begin{equation}\label{bfgs_proof_10}
(2\alpha_3\delta+\alpha_4)\frac{\eps}{1-r}\leq \delta
 \quad \text{and}\quad
 \gamma(1+r)[\tilde{L}\eps+2\eta\delta] \leq r.
\end{equation}
Based on the assumption that $\|\bbB_i^0-\nabla^2 f_i(\bbx^*)\|_\bbM\leq \delta$ we can derive the upper bound $\|\bbB_i^0-\nabla^2 f_i(\bbx^*)\|\leq \eta \delta$. This observation along with the inequality $\|\nabla^2 f_i(\bbx^*)\|\leq L$ implies that  $\|\bbB_i^0\|\leq \eta\delta+L$. Therefore, we obtain $\|(1/n)\sum_{i=1}^n\bbB_i^0\|\leq \eta\delta+L$. The second inequality in \eqref{bfgs_proof_10} implies that $2\gamma(1+r)\eta\delta\leq r$. Based on this observation and the inequalities $\|\bbB_i^0-\nabla^2 f_i(\bbx^*)\|\leq \eta \delta<2 \eta \delta$ and  $\gamma\geq\|\nabla^2 f_i(\bbx^*)^{-1}\|$, we obtain from Banach Lemma that $\|(\bbB_i^0)^{-1}\|\leq (1+r)\gamma$. Following the same argument for the matrix $((1/n)\sum_{i=1}^n \bbB_i^0)^{-1}$ with the inequalities $\|(1/n)\sum_{i=1}^n \bbB_i^0-(1/n)\sum_{i=1}^n \nabla^2 f_i(\bbx^*)\| 
\leq (1/n)\sum_{i=1}^n\| \bbB_i^0- \nabla^2 f_i(\bbx^*)\| \leq\eta \delta$ and  $\|\nabla^2 f(\bbx^*)^{-1}\|\leq \gamma$ we obtain that
\begin{align}\label{bfgs_proof_20}
\left\|\left(\frac{1}{n}\sum_{i=1}^n \bbB_i^0\right)^{-1}\right\|\leq (1+r)\gamma.
\end{align}
This upper bound in conjunction with the result in \eqref{lemma_claim_imp20} yields 
\begin{align}\label{bfgs_proof_30}
\|\bbx^{1}-\bbx^*\|& \leq   (1+r)\gamma \left[  \frac{\tilde{L}}{n} \sum_{i=1}^n  \left\|\bbz_i^0-\bbx^* \right\|^2+  \frac{1}{n} \sum_{i=1}^n  \left\| \left[\bbB_i^0-\nabla^2 f_{i}(\bbx^*)\right] \left(\bbz_i^0-\bbx^*\right) \right\|  \right]\nonumber\\
&=  (1+r)\gamma \left[  \tilde{L}  \left\|\bbx^{0}-\bbx^* \right\|^2\!\!+\!  \frac{1}{n} \sum_{i=1}^{n}  \left\| \left[\bbB_i^0-\nabla^2 f_{i}(\bbx^*)\right] \left(\bbx^{0}-\bbx^*\right) \right\|  \right].
\end{align}
Considering the assumptions that $\|\bbx^0-\bbx^*\|\leq \eps$ and $\|\bbB_i^0-\nabla^2 f_{i}(\bbx^*)\|\leq \eta \delta<2\eta \delta$ we can write
\begin{align}\label{bfgs_proof_40}
\|\bbx^{1}-\bbx^*\|&\leq (1+r)\gamma [\tilde{L}\eps+2\eta \delta] \|\bbx^{0}-\bbx^*\|
\nonumber\\
&\leq r\|\bbx^{0}-\bbx^*\|,
\end{align}
where the second inequality follows from the second condition in \eqref{bfgs_proof_10}. Without loss of generality, assume that $i_0=1$. Then, based on the result in \eqref{lemma_M_2_claim} we obtain
\begin{align}\label{bfgs_proof_50}
\left\|\bbB_1^{1}- \nabla^2 f_1(\bbx^*)\right\|_\bbM
& \leq \left[ (1- \alpha {\theta_1^0}^2)^{1/2}
	+\alpha_3 \sigma_1^0 \right]\left\|\bbB_1^{0}- \nabla^2 f_1(\bbx^*)\right\|_\bbM  + \alpha_4 \sigma_1^0\nonumber\\
&\leq (1+\alpha_3\eps) \delta +\alpha_4\eps\nonumber\\
&\leq \delta+ 2\alpha_3 \eps \delta +\alpha_4\eps \leq 2\delta. 
\end{align}
We proceed to the next iteration which leads to the inequality 
\begin{align}\label{bfgs_proof_60}
\|\bbx^{2}-\bbx^*\|& \leq   (1+r)\gamma \left[  \frac{\tilde{L}}{n} \sum_{i=1}^n  \left\|\bbz_i^t-\bbx^* \right\|^2+  \frac{1}{n} \sum_{i=1}^n  \left\| \left[\bbB_i^t-\nabla^2 f_{i}(\bbx^*)\right] \left(\bbz_i^t-\bbx^*\right) \right\|  \right]\nonumber\\
&\leq (1+r)\gamma\left[ {\tilde{L}} \eps +2\eta\delta \right] 
\left(  \frac{n-1}{n}\|\bbx^0-\bbx^*\|+\frac{1}{n}\|\bbx^1-\bbx^*\|\right)\nonumber\\
&\leq r
\left(  \frac{n-1}{n}\|\bbx^0-\bbx^*\|+\frac{1}{n}\|\bbx^1-\bbx^*\|\right)
\nonumber\\
&\leq r\|\bbx^0-\bbx^*\|.
\end{align}
And since the updated index is $i_1=2$ we obtain
\begin{align}\label{bfgs_proof_70}
\left\|\bbB_2^{2}- \nabla^2 f_2(\bbx^*)\right\|_\bbM
& \leq \left[ (1- \alpha {\theta_2^0}^2)^{1/2}
	+\alpha_3 \sigma_2^0 \right]\left\|\bbB_2^{0}- \nabla^2 f_2(\bbx^*)\right\|_\bbM  + \alpha_4 \sigma_2^0\nonumber\\
&\leq (1+\alpha_3\eps) \delta +\alpha_4\eps\nonumber\\
&\leq \delta+ 2\alpha_3 \eps \delta +\alpha_4\eps \leq 2\delta.
\end{align}
With the same argument we can show that all $\left\|\bbB_t^{t}- \nabla^2 f_t(\bbx^*)\right\|_\bbM\leq 2\delta$ and $\|\bbx^t-\bbx^*\|\leq \eps $, for all iterates $t=1,\dots, n$. Moreover, we have $\|\bbx^t-\bbx^*\|\leq r \|\bbx^0-\bbx^*\|$ for $t=1,\dots,n$. 

Now we use the results for iterates $t=1,\dots,n$ as the base of our induction argument. To be more precise, let's assume that for iterates $t=jn+1,jn+2,\dots, jn+n$ we know that the residuals are bounded above by $\|\bbx^t-\bbx^*\|\leq r^{j+1}\|\bbx^0-\bbx^*\|$ and the Hessian approximation matrices $\bbB_i^t$ satisfy the inequalities $\| \bbB_i^t- \nabla^2 f_i(\bbx^*)\|\leq 2\eta \delta$. Our goal is to show that for iterates $t=(j+1)n+1,(j+1)n+2,\dots, (j+1)n+n$ the inequalities $\|\bbx^t-\bbx^*\|\leq r^{j+2}\|\bbx^0-\bbx^*\|$ and $\| \bbB_i^t- \nabla^2 f_i(\bbx^*)\|\leq 2\eta \delta$ hold. 

Based on the inequalities $\| \bbB_i^t- \nabla^2 f_i(\bbx^*)\|\leq 2\eta \delta$ and $\|\nabla^2 f_i(\bbx^*)^{-1}\|\leq \gamma$ we can show that for all $t=jn+1,jn+2,\dots, jn+n$ we have
\begin{align}\label{bfgs_proof_80}
\left\|\left(\frac{1}{n}\sum_{i=1}^n \bbB_i^t\right)^{-1}\right\|\leq (1+r)\gamma.
\end{align}
Using \eqref{bfgs_proof_80} and the inequality in \eqref{lemma_M_2_claim} for the iterate $t=(j+1)n+1$, we obtain 
\begin{align}\label{bfgs_proof_90}
\|\bbx^{(j+1)n+1}-\bbx^*\|& \leq   (1+r)\gamma  \frac{\tilde{L}}{n} \sum_{i=1}^n  \left\|\bbz_i^{(j+1)n}-\bbx^* \right\|^2
\nonumber\\
&\quad 
+(1+r)\gamma   \frac{1}{n} \sum_{i=1}^n  \left\| \left[\bbB_i^{(j+1)n}-\nabla^2 f_{i}(\bbx^*)\right] \left(\bbz_i^{(j+1)n}-\bbx^*\right) \right\|.
\end{align}
Since the variables are updated in a cyclic fashion the set of variables $\{\bbz_i^{(j+1)n}\}_{i=1}^{i=n}$ is equal to the set $\{\bbx^{(j+1)n-i}\}_{i=0}^{i=n-1}$. By considering this relation and replacing the norms  $\| [\bbB_i^{(j+1)n}-\nabla^2 f_{i}(\bbx^*)](\bbz_i^{(j+1)n}-\bbx^*) \|$ by their upper bounds $2\eta\delta\|\bbz_i^{(j+1)n}-\bbx^*\|$ we can simplify the right hand side of \eqref{bfgs_proof_90} as
\begin{align}\label{bfgs_proof_100}
\|\bbx^{(j+1)n+1}-\bbx^*\| \leq   (1+r)\gamma \left[  \frac{\tilde{L}}{n} \sum_{i=1}^n  \left\|\bbx^{jn+i}-\bbx^* \right\|^2+   \frac{2\eta\delta}{n} \sum_{i=1}^n  \left\| \bbx^{jn+i}-\bbx^* \right\|\right].
\end{align}
Since $\|\bbx^{jn+i}-\bbx^* \|\leq \eps$ for all $j=1,\dots,n$, we obtain 
\begin{align}\label{bfgs_proof_110}
\|\bbx^{(j+1)n+1}-\bbx^*\| \leq   (1+r)\gamma\left[ {\tilde{L}} \eps +2\eta\delta \right] \left(\frac{1}{n}\sum_{i=1}^n  \left\| \bbx^{jn+i}-\bbx^* \right\|\right).
\end{align}
According to the second inequality in \eqref{bfgs_proof_10} and the assumption that for iterates $t=jn+1,jn+2,\dots, jn+n$ we know that $\|\bbx^t-\bbx^*\|\leq r^{j+1}\|\bbx^0-\bbx^*\|$, we can replace the right hand side of \eqref{bfgs_proof_110} by the following upper bound 
\begin{align}\label{bfgs_proof_120}
\|\bbx^{(j+1)n+1}-\bbx^*\|\leq r^{j+2}\|\bbx^0-\bbx^*\|.
\end{align}
Now we show that the updated Hessian approximation $ \bbB_{i_t}^{(j+1)n+1}$ for $t=(j+1)n+1$ satisfies the inequality $\|\bbB_{i_t}^{(j+1)n+1}- \nabla^2 f_{i_t}(\bbx^*)\|_\bbM\leq 2\delta$. According to the result in \eqref{lemma_M_2_claim}, we can write 
\begin{align}\label{bfgs_proof_130}
&\left\|\bbB_{i_t}^{(j+1)n+1}- \nabla^2 f_{i_t}(\bbx^*)\right\|_\bbM
-\left\|\bbB_{i_t}^{jn+1}- \nabla^2 f_{i_t}(\bbx^*)\right\|_\bbM
\nonumber\\ 
&\qquad \qquad \leq \alpha_3\sigma_{i_t}^{jn+1} \left\|\bbB_{i_t}^{jn+1}- \nabla^2 f_{i_t}(\bbx^*)\right\|_\bbM
 + \alpha_4 \sigma_{i_t}^{jn+1}.
\end{align}
Now observe that $\sigma_{i_t}^{jn+1}=\max\{\|\bbx^{(j+1)n+1}-\bbx^*\|,\|\bbx^{jn+1}-\bbx^*\|\}$ is bounded above by $r^{j+1}\|\bbx^0-\bbx^*\|$. Applying this substitution into \eqref{bfgs_proof_130} and considering the conditions $\|\bbB_{i_t}^{jn+1}- \nabla^2 f_{i_t}(\bbx^*)\|_\bbM\leq 2\delta$ and $\|\bbx^0-\bbx^*\|\leq \eps$ lead to the inequality 
\begin{align}\label{bfgs_proof_140}
\left\|\bbB_{i_t}^{(j+1)n+1}- \nabla^2 f_{i_t}(\bbx^*)\right\|_\bbM
-\left\|\bbB_{i_t}^{jn+1}- \nabla^2 f_{i_t}(\bbx^*)\right\|_\bbM
 \leq r^{j+1}\eps( 2\delta\alpha_3 + \alpha_4 ).
\end{align}
By writing the expression in \eqref{bfgs_proof_140} for previous iterations and using a recursive logic we obtain that 
\begin{align}\label{bfgs_proof_150}
\left\|\bbB_{i_t}^{(j+1)n+1}- \nabla^2 f_{i_t}(\bbx^*)\right\|_\bbM
-\left\|\bbB_{i_t}^{0}- \nabla^2 f_{i_t}(\bbx^*)\right\|_\bbM
 \leq \eps( 2\delta\alpha_3 + \alpha_4 )\frac{1}{1-r}.
\end{align}
Based on the first inequality in \eqref{bfgs_proof_10}, the right hand side of \eqref{bfgs_proof_150} is bounded above by $\delta$. Moreover, the norm $\|\bbB_{i_t}^{0}- \nabla^2 f_{i_t}(\bbx^*)\|_\bbM$ is also upper bounded by $\delta$. These two bounds imply that 
\begin{align}\label{bfgs_proof_160}
\left\|\bbB_{i_t}^{(j+1)n+1}- \nabla^2 f_{i_t}(\bbx^*)\right\|_\bbM\leq 2\delta,
\end{align}
and consequently $\|\bbB_{i_t}^{(j+1)n+1}- \nabla^2 f_{i_t}(\bbx^*)\|\leq 2\eta \delta$. By following the steps from \eqref{bfgs_proof_90} to \eqref{bfgs_proof_160}, we can show for all iterates $t=(j+1)n+1,(j+1)n+2,\dots, (j+1)n+n$ the inequalities $\|\bbx^t-\bbx^*\|\leq r^{j+2}\|\bbx^0-\bbx^*\|$ and $\| \bbB_i^t- \nabla^2 f_i(\bbx^*)\|\leq 2\eta \delta$ hold. The induction proof is complete and \eqref{loc_lin_convg55} holds. Moreover, the inequality $\| \bbB_i^t- \nabla^2 f_i(\bbx^*)\|\leq 2\eta \delta$ holds for all $i$ and steps $t$. Hence, the norms $\|\bbB_i^t\|$ and $\|({\bbB_i^t})^{-1}\|$, and consequently $\|(1/n)\sum_{i=1}^n\bbB_i^t\|$ and $\|((1/n)\sum_{i=1}^n\bbB_i^t)^{-1}\|$ are uniformly bounded.

\section{Proof of Proposition \ref{prop_DM_condition}}\label{apx_prop_DM_condition}
According to the result in Lemma \ref{theorem_local_lin_conv}, we can show that the sequence  of errors $\sigma_i^t=\max\{\|\bbz_i^{t+1}-\bbx^*\|,\|\bbz_i^t-\bbx^*\|\}$ is summable for all $i$. To do so, consider the sum of the sequence $\sigma_i^t$ which is upper bounded by
\begin{align}\label{proof_theorem_eq_10}
\sum_{t=0}^\infty\sigma_i^t
&= \sum_{t=0}^\infty  \max\{\|\bbz_i^{t+1}-\bbx^*\|,\|\bbz_i^t-\bbx^*\|\} \leq \sum_{t=0}^\infty \|\bbz_i^{t+1}-\bbx^*\|+  \sum_{t=0}^\infty \|\bbz_i^t-\bbx^*\|
\end{align}
Note that the last time that the index $i$ is chosen before time $t$ should be in the set $\{t-1, \dots, t-n\}$. This observation in association with the result in \eqref{loc_lin_convg55} implies that 
\begin{align}\label{proof_theorem_eq_20}
\sum_{t=0}^\infty\sigma_i^t \leq 2\sum_{t=0}^\infty r^{[\frac{t-n-1}{n}]+1} \|\bbx^0-\bbx^*\| = 2\sum_{t=0}^\infty r^{[\frac{t-1}{n}]} \|\bbx^0-\bbx^*\|
\end{align}
Simplifying the sum in the right hand side of \eqref{proof_theorem_eq_20} yields 
\begin{align}\label{proof_theorem_eq_30}
\sum_{t=0}^\infty\sigma_i^t 
 \leq \frac{2\|\bbx^0-\bbx^*\|}{r}+ 2n\|\bbx^0-\bbx^*\|\sum_{t=0}^\infty r^{t} <\infty .
\end{align}
Thus, the sequence $\sigma_i^t$ is summable for all $i=1,\dots,n$. To complete the proof we use the following result from Lemma 3.3 in \cite{Dennis}.

\begin{lemma}\label{lemma_dennis}
Let $\{\phi^t\}$ and $\{\delta^t\}$ be sequences of nonnegative numbers such that 
\begin{equation}
\phi^{t+1} \leq (1+\delta^t)\phi^t +\delta^t 
\quad 
\text{and} 
\quad\sum_{k=1}^\infty \delta^t<\infty.
\end{equation}
Then, the sequence $\{\phi^t\}$ converges.
\end{lemma}

Considering the results in Lemmata \ref{lemma_M_2} and \ref{lemma_dennis}, and the fact that $\sigma_i^t$ is summable as shown in\eqref{proof_theorem_eq_30}, we obtain that the sequence $\left\|\bbB_i^t- \nabla^2 f_i(\bbx^*)\right\|_\bbM$ for $\bbM:= \nabla^2 f_i(\bbx^*)^{-1/2}$ is convergent and the following limit exists 
\begin{equation}\label{proof_theorem_eq_40}
\lim_{k\to\infty} \| \nabla^2 f_i(\bbx^*)^{-1/2}\  \bbB_i^t\ \nabla^2 f_i(\bbx^*)^{-1/2} \ -\ \bbI\|_\bbF=l
\end{equation}
where $l$ is a nonnegative constant. Moreover, following the proof of Theorem 3.4 in \cite{Dennis} we can show that
\begin{align}\label{proof_theorem_eq_50}
\alpha (\theta_i^t)^2 \|\bbB_i^t- \nabla^2 f_i(\bbx^*)\|_\bbM
& \leq  \|\bbB_i^t- \nabla^2 f_i(\bbx^*)\|_\bbM -  \|\bbB_i^{t+1}- \nabla^2 f_i(\bbx^*)\|_\bbM \nonumber\\
&\quad +\sigma_i^t (\alpha_3  \|\bbB_i^t- \nabla^2 f_i(\bbx^*)\|_\bbM +\alpha_4), 
\end{align}
and, therefore, summing both sides implies, 
\begin{align}\label{proof_theorem_eq_60}
\sum_{t=0}^\infty (\theta_i^t)^2 \|\bbB_i^t- \nabla^2 f_i(\bbx^*)\|_\bbM < \infty
\end{align}
Replacing $\theta_i^t$ in \eqref{proof_theorem_eq_60} by its definition in \eqref{def_theta} results in 
\begin{align}\label{proof_theorem_eq_70}
\sum_{t=0}^\infty\frac{\|\bbM(\bbB_i^t- \nabla^2 f_i(\bbx^*))\bbs_i^t\|^2}{\|\bbB_i^t- \nabla^2 f_i(\bbx^*)\|_\bbM\|\bbM^{-1}\bbs_i^t\|^2} < \infty
\end{align}
Since the norm $\|\bbB_i^t- \nabla^2 f_i(\bbx^*)\|_\bbM$ is upper bounded and the eigenvalues of the matrix $\bbM= \nabla^2 f_i(\bbx^*)^{-1/2} $ are uniformly lower and upper bounded, we conclude from the result in \eqref{proof_theorem_eq_70} that 
\begin{align}\label{proof_theorem_eq_80}
\lim_{t\to\infty} \frac{\|(\bbB_i^t- \nabla^2 f_i(\bbx^*))\bbs_i^t\|^2}{\|\bbs_i^t\|^2}=0,
\end{align}
which yields the claim in \eqref{instan_dennis_more}.

\section{Proof of Lemma  \ref{lemma_incr_bfgs_lim_final}}\label{apx_lemma_incr_bfgs_lim_final} Consider the sets of variable variations $\ccalS_1=\{\bbs_i^{t+n\tau}\}_{\tau=0}^{\tau=T}$ and $\ccalS_2=\{\bbs_i^{t+n\tau}\}_{\tau=0}^{\tau=\infty}$. It is trivial to show that $\bbz_i^{t}-\bbx^*$ is in the span of the set $\ccalS_2$, since the sequences of variables $\bbx^t$ and $\bbz_i^t$ converge to $\bbx^*$ and we can write $\bbx^*-\bbz_i^{t}= \sum_{\tau=0}^{\infty} \bbs_i^{t+n\tau}$. We proceed to show that the vector $\bbz_i^t-\bbx^*$ is also in the span of the set $\ccalS_1$ when $T$ is sufficiently large. To do so, we use a contradiction argument. Let's assume that the vector $\bbz_i^t-\bbx^*$ does not lie in the span of the set $\ccalS_1$, and,  therefore, it can be decomposed as the sum of two non-zero vectors given by 
\begin{equation}
\bbz_i^t-\bbx^*= \bbv^{t}_{\parallel} +\bbv^{t}_{\perp},
\end{equation}
where $\bbv^{t}_{\parallel}$ lies in the span of $\ccalS_1$ and $\bbv^{t}_{\perp}$ is orthogonal to the span of $\ccalS_1$. Since we assume that $\bbz_i^t-\bbx^*$ does not lie in the span of $\ccalS_1$, we obtain that $\bbz_i^{t+nT}-\bbx^*$ also does not lie in this span, since $\bbz_i^{t+nT}-\bbx^*$ can be written as the sum $\bbz_i^{t+nT}-\bbx^*= \bbz_i^t-\bbx^*+\sum_{\tau=0}^T\bbs_i^{t+n\tau}$. These observations imply that we can also decompose the vector $\bbz_i^{t+nT}-\bbx^*$ as
\begin{equation}
\bbz_i^{t+nT}-\bbx^*= \bbv^{t+nT}_{\parallel} +\bbv^{t+nT}_{\perp},
\end{equation}
where $\bbv^{t+nT}_{\parallel}$ lies in the span of $\ccalS_1$ and $\bbv^{t+nT}_{\perp}$ is orthogonal to the span of $\ccalS_1$. Moreover, we obtain that $\bbv^{t+nT}_{\perp}$ is equal to $\bbv^{t}_{\perp}$, i.e., 
\begin{equation}
\bbv^{t+nT}_{\perp}=\bbv^{t}_{\perp}.
\end{equation}
This is true since $\bbz_i^{t+nT}-\bbx^*$ can be written as the sum of $\bbz_i^{t}-\bbx^*$ and a group of vectors that lie in the span of $\ccalS_1$. We assume that the norm $\|\bbv^{t+nT}_{\perp}\|=\|\bbv^{t}_{\perp}\|=\eps$ where $\eps>0$ is a strictly positive constant. According to the linear convergence of the sequence $\|\bbx^t-\bbx^*\|$ in Lemma \ref{theorem_local_lin_conv} we know that 
\begin{equation}
\|\bbz_i^{t+nT}-\bbx^*\|\leq r^{[\frac{t+nT-1}{n}]+1}  \|\bbx^0-\bbx^*\|\leq r^{T}\|\bbx^0-\bbx^*\|
\end{equation}
If we pick large enough $T$ such that $r^{T}\|\bbx^0-\bbx^*\|<\eps$, then we obtain $\|\bbz_i^{t+nT}-\bbx^*\|<\eps$ which contradicts the assumption $\|\bbv^{t}_{\perp}\|=\eps$. Thus, we obtain that the vector $\bbz_i^{t}-\bbx^*$ is also in the span of set $\ccalS_1$. 

Since the vector $\bbz_i^{t}-\bbx^*$ is in the span of $\ccalS_1$, we can write the normalized vector $(\bbz_i^{t}-\bbx^*)/\|\bbz_i^{t}-\bbx^*\|$ as a linear combination of the set of normalized vectors $\{\bbs_i^{t+n\tau}/\|\bbs_i^{t+n\tau}\|\}_{\tau=0}^{\tau=T}$. This property allows to write 
\begin{align}\label{opop_10}
\lim_{t \rightarrow \infty} \frac{ \| (\bbB_i^t - \nabla^2 f_i(\bbx^*)) (\bbz_i^{t} - \bbx^*) \|}{\| \bbz_i^{t} -\bbx^*\|} 
&=
\lim_{t \rightarrow \infty} \left\| (\bbB_i^t - \nabla^2 f_i(\bbx^*))  \frac{(\bbz_i^{t}- \bbx^*) }{\| \bbz_i^{t} -\bbx^*\|} \right\|
\nonumber\\
&
=
\lim_{t \rightarrow \infty} \left\| (\bbB_i^t - \nabla^2 f_i(\bbx^*))  \sum_{\tau=0}^T a_\tau \frac{ \bbs_i^{t+n\tau} }{\|\bbs_i^{t+n\tau}\|}\right\|,
\end{align}
where $a_\tau$ is coefficient of the vector $\bbs_i^{t+n\tau}$ when we write $(\bbz_i^{t} - \bbx^*)/\|\bbz_i^{t}- \bbx^*\|$ as the linear combination of the normalized vectors $\{\bbs_i^{t+n\tau}/\|\bbs_i^{t+n\tau}\|\}_{\tau=0}^{\tau=T}$. Now since the index of the difference $\bbB_i^t - \nabla^2 f_i(\bbx^*)$ does not match with the descent directions $\bbs_i^t+n\tau$. We add and subtract the term $\bbB_i^{t+n\tau} $ to the expression $\bbB_i^t - \nabla^2 f_i(\bbx^*)$ and use the triangle inequality to write 
\begin{align}\label{opop_20}
&\lim_{t \rightarrow \infty} \frac{ \| (\bbB_i^t - \nabla^2 f_i(\bbx^*)) (\bbz_i^{t}- \bbx^*) \|}{\| \bbz_i^{t} -\bbx^*\|} \nonumber\\
&\leq
\lim_{t \rightarrow \infty} \left\|   \sum_{\tau=0}^T a_\tau \frac{(\bbB_i^{t+n\tau}-\nabla^2 f_i(\bbx^*)) \bbs_i^{t+n\tau} }{\|\bbs_i^{t+n\tau}\|}\right\|+\left\|   \sum_{\tau=0}^T a_\tau \frac{(\bbB_i^t - \bbB_i^{t+n\tau})\bbs_i^{t+n\tau} }{\|\bbs_i^{t+n\tau}\|}\right\|.
\end{align}
We first simplify the first limit in the right hand side of \eqref{opop_20}. Using the Cauchy-Schwarz inequality and the result in Proposition \ref{prop_DM_condition} we can write 
\begin{align}\label{opop_30}
\lim_{t \rightarrow \infty} \left\|   \sum_{\tau=0}^T a_\tau \frac{(\bbB_i^{t+n\tau}\!-\!\nabla^2 f_i(\bbx^*)) \bbs_i^{t+n\tau} }{\|\bbs_i^{t+n\tau}\|}\right\|
&
\leq \lim_{t \rightarrow \infty}   \sum_{\tau=0}^T a_\tau\left\|  \frac{(\bbB_i^{t+n\tau}\!-\!\nabla^2 f_i(\bbx^*)) \bbs_i^{t+n\tau} }{\|\bbs_i^{t+n\tau}\|}\right\|
\nonumber\\
&=   \sum_{\tau=0}^T a_\tau \! \lim_{t \rightarrow \infty} \left\|  \frac{(\bbB_i^{t+n\tau}\!-\!\nabla^2 f_i(\bbx^*)) \bbs_i^{t+n\tau} }{\|\bbs_i^{t+n\tau}\|}\right\|
=
0.
\end{align}
Based on the results in \eqref{opop_20} and \eqref{opop_30}, to prove the claim in \eqref{eq_lemma_incr_bfgs_lim_final} it remains to show
\begin{align}\label{opop_40}
\lim_{t \rightarrow \infty} \left\|   \sum_{\tau=0}^T a_\tau \frac{(\bbB_i^t - \bbB_i^{t+n\tau})\bbs_i^{t+n\tau} }{\|\bbs_i^{t+n\tau}\|}\right\| =0.
\end{align}

To reach this goal, we first study the limit of the difference between two consecutive update Hessian approximation matrices $\lim_{t \rightarrow \infty} \|\bbB_i^t-\bbB_i^{t+n}\|$. Note that if we set $\bbA=\bbB_i^t $ in \eqref{B_cons_rela}, we obtain that 
\begin{align}\label{opop_50}
\|\bbB_i^{t+n}-\bbB_i^t\|_\bbM \leq \alpha_2 \frac{\|\bby_i^t-\bbB_i^t\bbs_i^t\|}{\|\bbM^{-1}\bbs_i^t\|}.
\end{align} 
where $\bbM=(\nabla^2 f_i(\bbx^*))^{-1/2}$. By adding and subtracting the term $\nabla^2 f_i(\bbx^*)\bbs_i^t$ and using the result in \eqref{instan_dennis_more}, we can show that the difference $\|\bbB_i^{t+n}-\bbB_i^t\|_\bbM$ approaches zero asymptotically. In particular, 
\begin{align}\label{opop_60}
\lim_{t \rightarrow \infty} \|\bbB_i^{t+n}-\bbB_i^t\|_\bbM 
&\leq \alpha_2 \lim_{t \rightarrow \infty} \frac{\|\bby_i^t-\bbB_i^t\bbs_i^t\|}{\|\bbM^{-1}\bbs_i^t\|}
\nonumber\\
&\leq  \alpha_2 \lim_{t \rightarrow \infty} \frac{\|\bby_i^t- \nabla^2 f_i(\bbx^*)\bbs_i^t\|}{\|\bbM^{-1}\bbs_i^t\|}+
\alpha_2 \lim_{t \rightarrow \infty} \frac{\|(\nabla^2 f_i(\bbx^*)-\bbB_i^t)\bbs_i^t\|}{\|\bbM^{-1}\bbs_i^t\|}.
\end{align} 
Since $\|\bby_i^t- \nabla^2 f_i(\bbx^*)\bbs_i^t\|$ is bounded above by $\tilde{L}\|\bbs_i^t\|\max\{\|\bbz_i^{t}-\bbx^*\|,\|\bbz_i^{t+1}-\bbx^*\|\}$ and the eigenvalues of the matrix $\bbM$ are uniformly bounded we obtain that the first limit in the right hand side of \eqref{opop_60} converges to zero. Further, the result in \eqref{instan_dennis_more} shows that the second limit in the right hand side of \eqref{opop_60} also converges to zero. Therefore,  
\begin{align}\label{opop_70}
\lim_{t \rightarrow \infty} \|\bbB_i^{t+n}-\bbB_i^t\|_\bbM =0.
\end{align}
Following the same argument we can show that for any two consecutive Hessian approximation matrices the difference approaches zero asymptotically. Thus, we obtain 
\begin{align}\label{opop_80}
\lim_{t \rightarrow \infty} \left\|\bbB_i^t - \bbB_i^{t+n\tau}\right\|_\bbM 
&\leq 
\lim_{t \rightarrow \infty} \left\|\sum_{u=0}^{\tau-1}\left(\bbB_i^{t+nu} - \bbB_i^{t+n(u+1)}\right)\right\|_\bbM 
\nonumber\\
&
\leq 
\sum_{u=0}^{\tau-1} \lim_{t \rightarrow \infty} \left\|\bbB_i^{t+nu} - \bbB_i^{t+n(u+1)}\right\|_\bbM =0.
\end{align}
Observing the result in \eqref{opop_80} we can show that 
\begin{align}\label{opop_90}
\lim_{t \rightarrow \infty} \left\|   \sum_{\tau=0}^T a_\tau \frac{(\bbB_i^t - \bbB_i^{t+n\tau})\bbs_i^{t+n\tau} }{\|\bbs_i^{t+n\tau}\|}\right\| 
&\leq 
\sum_{\tau=0}^T a_\tau\lim_{t \rightarrow \infty} \left\|    \frac{(\bbB_i^t - \bbB_i^{t+n\tau})\bbs_i^{t+n\tau} }{\|\bbs_i^{t+n\tau}\|}\right\| 
\nonumber\\
&\leq 
\sum_{\tau=0}^T a_\tau\lim_{t \rightarrow \infty} \left\| \bbB_i^t - \bbB_i^{t+n\tau}\right\| =0. 
\end{align}
Therefore, the result in \eqref{opop_40} holds. The claim in \eqref{eq_lemma_incr_bfgs_lim_final} follows by combining the results in \eqref{opop_20}, \eqref{opop_30}, and \eqref{opop_40}.

\section{Proof of Theorem  \ref{sup_linear_thm}}\label{apx_sup_linear_thm} The result in Lemma \ref{lemma_cyclic_property} implies
\begin{align}\label{lemma_claim_100}
&\|\bbx^{t+1}-\bbx^*\|
\leq
 \frac{\tilde{L}\Gamma^{t}}{n} \sum_{i=1}^n  \left\|\bbz_i^t-\bbx^* \right\|^2
 +
 \frac{\Gamma^{t}}{n} \sum_{i=1}^n  \left\| \left(\bbB_i^t-\nabla^2 f_{i}(\bbx^*)\right) \left(\bbz_i^t-\bbx^*\right) \right\|.
\end{align}
Divide both sides of \eqref{lemma_claim_100} by $(1/n) \sum_{i=1}^n  \left\|\bbz_i^t-\bbx^* \right\|$ to obtain 
\begin{equation}\label{lemma_claim_200}
\frac{\|\bbx^{t+1}-\bbx^*\|}{ \frac{1}{n}\sum_{i=1}^n  \left\|\bbz_i^t-\bbx^* \right\|}
\leq{\tilde{L}\Gamma^{t}} \sum_{i=1}^n \frac{  \left\|\bbz_i^t-\bbx^* \right\|^2}{ \sum_{i=1}^n  \left\|\bbz_i^t-\bbx^* \right\|}
 +
 {\Gamma^{t}}\sum_{i=1}^n \frac{  \left\| \left(\bbB_i^t-\nabla^2 f_{i}(\bbx^*)\right) \left(\bbz_i^t-\bbx^*\right) \right\|}{ \sum_{i=1}^n  \left\|\bbz_i^t-\bbx^* \right\|}
\end{equation}
Since the error $\|\bbz_i^t-\bbx^*\|$ is a lower bound for the sum of errors $\sum_{i=1}^n  \|\bbz_i^t-\bbx^* \|$, we can replace $\|\bbz_i^t-\bbx^*\|$ for $\sum_{i=1}^n  \|\bbz_i^t-\bbx^* \|$ into \eqref{lemma_claim_200} which implies
\begin{align}\label{lemma_claim_300}
 \frac{\|\bbx^{t+1}-\bbx^*\|}{ \frac{1}{n}\sum_{i=1}^n  \left\|\bbz_i^t-\bbx^* \right\|}
& \leq
 {\tilde{L}\Gamma^{t}} \sum_{i=1}^n \frac{  \left\|\bbz_i^t-\bbx^* \right\|^2}{   \left\|\bbz_i^t-\bbx^* \right\|}
 +
 {\Gamma^{t}} \sum_{i=1}^n \frac{  \left\| \left(\bbB_i^t-\nabla^2 f_{i}(\bbx^*)\right) \left(\bbz_i^t-\bbx^*\right) \right\|}{  \left\|\bbz_i^t-\bbx^* \right\|}
  \nonumber\\
  &
  =
  {\tilde{L}\Gamma^{t}} \sum_{i=1}^n {  \left\|\bbz_i^t-\bbx^* \right\|}
   +
 {\Gamma^{t}} \sum_{i=1}^n \frac{  \left\| \left(\bbB_i^t-\nabla^2 f_{i}(\bbx^*)\right) \left(\bbz_i^t-\bbx^*\right) \right\|}{  \left\|\bbz_i^t-\bbx^* \right\|}.
\end{align}
Since $\Gamma^t$ is bounded above, computing the limit of both sides in \eqref{lemma_claim_300} yields
\begin{align}\label{lemma_claim_200000}
\lim_{t\to\infty}\frac{\|\bbx^{t+1}-\bbx^*\|}{ \frac{1}{n}\sum_{i=1}^n  \left\|\bbz_i^t-\bbx^* \right\|}  = 0.
\end{align}
The result in \eqref{lemma_claim_200000} in association with the simplification for the sum $\sum_{i=1}^n  \left\|\bbz_i^t-\bbx^* \right\|=\sum_{i=0}^{n-1}  \left\|\bbx^{t-i}-\bbx^* \right\|$ leads to the claim in \eqref{final_claim}.

\section{Proof of Theorem  \ref{sup_linear_thm_2}}\label{apx_sup_linear_thm_2} Consider the definition of the sequence  $\tbx^t=\argmax_{u\in\{tn,\dots,tn+n-1\}}\{\|\bbx^u-\bbx^*\|\}$ which is a subsequence of the sequence $\{\bbx^t\}_{t=0}^{\infty}$. Our goal is to show this subsequence converges superlinearly to $\bbx^*$, i.e., 
$\lim_{t\to\infty} \frac{\|\tbx^{t+1}-\bbx^*\|}{\|\tbx^{t}-\bbx^*\|} =0.$
To do so, first note that the result in Theorem \ref{sup_linear_thm} implies that 
\begin{align}\label{sup_upp_bound_proof_200}
\lim_{t\to\infty}\ \frac{\|\bbx^t-\bbx^*\|}{\max\{\|\bbx^{t-1}-\bbx^*\| ,\dots, \|\bbx^{t-n}-\bbx^*\| \}}=0,
\end{align}
 which follows from the inequality $\max\{\|\bbx^{t-1}-\bbx^*\| ,\dots, \|\bbx^{t-n}-\bbx^*\| \}\geq (1/n)(\|\bbx^{t-1}-\bbx^*\|+\dots+\|\bbx^{t-n}-\bbx^*\|)$.  Based on the limit in \eqref{sup_upp_bound_proof_200}, there exists a large enough $t_0$ such that for all $t\geq t_0$ the following inequality holds,
\begin{align}\label{sup_upp_bound_proof_300}
\|\bbx^{t}-\bbx^*\| < \max\{\|\bbx^{t-1}-\bbx^*\| ,\dots, \|\bbx^{t-n}-\bbx^*\| \}.
\end{align}
Combining the inequality in \eqref{sup_upp_bound_proof_300} with the inequalities $\|\bbx^{t-i}-\bbx^*\| \leq \max\{\|\bbx^{t-1}-\bbx^*\| ,\dots, \|\bbx^{t-n}-\bbx^*\| \}$ for $i=1,\dots,n-1$ yields
\begin{align}\label{sup_upp_bound_proof_400}
 \max\{\|\bbx^{t}-\bbx^*\| ,\dots, \|\bbx^{t-n+1}-\bbx^*\| \}\leq \max\{\|\bbx^{t-1}-\bbx^*\| ,\dots, \|\bbx^{t-n}-\bbx^*\| \},
\end{align}
and consequently we can generalize this result to obtain
\begin{align}\label{sup_upp_bound_proof_500}
 \max\{\|\bbx^{t}\!-\!\bbx^*\| ,\dots, \|\bbx^{t-n+1}\!-\!\bbx^*\| \}\leq \max\{\|\bbx^{t-\tau}\!-
\!\bbx^*\| ,\dots, \|\bbx^{t-\tau-n+1}\!-\!\bbx^*\| \},
\end{align}
for any positive integer $\tau$ such that $t-\tau\geq t_0$. 

We use the result in \eqref{sup_upp_bound_proof_500} to build a superlinearly convergent subsequence of the residuals sequence $\|\bbx^t-\bbx^*\|$. If we define $\bbx^{tn+u_t^*}$ as the iterate that has the largest error among the iterates in the $t+1$-th pass, i.e.,
\begin{align}\label{sup_upp_bound_proof_600}
 \bbx^{tn+u_t^*}=\argmax_{u\in\{tn,\dots,tn+n-1\}}\{\|\bbx^u-\bbx^*\|\},
\end{align}
then it follows that $\tbx^t=\bbx^{tn+u_t^*}$, where $u_t^*\in \{0,1,\dots,n-1\}$.  Moreover, we obtain
\begin{align}\label{sup_upp_bound_proof_700}
 \frac{\|\tbx^{t+1}-\bbx^*\|}{\|\tbx^{t}-\bbx^*\|} 
&=
\frac{ \|\bbx^{(t+1)n+u_{t+1}^*}-\bbx^*\|}
{ \max\{\|\bbx^{tn}-\bbx^*\| ,\dots, \|\bbx^{tn+n-1}-\bbx^*\| \}}
\nonumber\\
& \leq 
\frac{ \|\bbx^{tn+n+u_{t+1}^*}-\bbx^*\|}
{ \max\{\|\bbx^{tn+u_{t+1}^*-1}-\bbx^*\| ,\dots, \|\bbx^{tn+n+u_{t+1}^*-1}-\bbx^*\| \}}.
\end{align}
The equality follows from the definition of the iterate $ \tbx^t$ and the definition in \eqref{sup_upp_bound_proof_600}, and the inequality holds because of the result in \eqref{sup_upp_bound_proof_500}. Considering the result in \eqref{sup_upp_bound_proof_200}, computing the limit of both sides leads to the conclusion that the sequence $\|\tbx^{t}-\bbx^*\|$ is superlinearly convergent. In other words, we obtain that the subsequence $\{\|\bbx^{tn+u_{t}^*}-\bbx^*\|\}_{t=0}^{t=\infty}$ superlinearly converges to zero.

Let's define the sequence $q^t$ such that $q^{kn}=\dots=q^{kn+n-1}=\|\tbx^k-\bbx^*\|$ for $k=0,1,2, \dots$, which means that the value of the sequence $q^t$ is fixed for each pass and is equal to the max error of the corresponding pass. Therefore, it is trivial to show that $q^t$ is always larger than or equal to $\|\bbx^t-\bbx^*\|$, i.e., $\|\bbx^t-\bbx^*\|\leq q^t$ for all $t\geq 0$. Now define the sequence $\zeta^t$ such that $\zeta^t=q^t$ for $t=0,\dots,n-1$, and for $t\geq n$
\begin{align}\label{sup_upp_bound_proof_800}
\zeta^{kn+i}= q^{kn-1}\left(\frac{q^{kn+n-1}}{q^{kn-1}}\right)^{\frac{i+1}{n}}, \qquad \for \quad i=0,\dots, n-1 , \quad k\geq 1.
\end{align}
According to this definition we can verify that $\zeta^t$ is an upper bound for the sequence $q^t$ and, consequently, an upper bound for the sequence of errors $\|\bbx^t-\bbx^*\|$. Based on the definition of the sequence $\zeta^t$ in \eqref{sup_upp_bound_proof_800}, the ratio $\zeta^{t+1}/\zeta^{t}$ is given by $(q^{\left\lfloor \frac{t+1}{n} \right\rfloor n+n-1}/q^{\left\lfloor \frac{t+1}{n} \right\rfloor n-1})^{1/n}$. This simplification in association with the definitions of the sequences $\|\bbx^t-\bbx^*\|$ and $\|\tbx^t-\bbx^*\|$ implies that 
\begin{align}\label{sup_upp_bound_proof_900}
\lim_{t\to\infty}\ \frac{\zeta^{t+1}}{\zeta^{t}}\ 
= \ \lim_{t\to\infty}\ \left(\frac{q^{\left\lfloor \frac{t+1}{n} \right\rfloor n+n-1}}{q^{\left\lfloor \frac{t+1}{n} \right\rfloor n-1}}\right)^{\frac{1}{n}}\
= \ \lim_{t\to\infty}\ \left(\frac{\|\tbx^{{\left\lfloor \frac{t+1}{n} \right\rfloor }}-\bbx^*\|}{\|\tbx^{{\left\lfloor \frac{t+1}{n} \right\rfloor-1 }}-\bbx^*\|}\right)^{\frac{1}{n}} \ = \ 0,
\end{align}
which leads to the claim in \eqref{sup_linear_thm_2_claim}.

\bibliographystyle{siamplain}
\bibliography{bmc_article,bmc_article2,bibliography,bibliography2}

\begin{thebibliography}{10}

\bibitem{blatt2007convergent}
{\sc D.~Blatt, A.~O. Hero, and H.~Gauchman}, {\em A convergent incremental
  gradient method with a constant step size}, SIAM Journal on Optimization, 18
  (2007), pp.~29--51.

\bibitem{Bottou}
{\sc L.~Bottou}, {\em Large-scale machine learning with stochastic gradient
  descent}, In Proceedings of COMPSTAT'2010,  (2010), pp.~177--186.

\bibitem{BottouCun}
{\sc L.~Bottou and Y.~L. Cun}, {\em On-line learning for very large datasets},
  in Applied Stochastic Models in Business and Industry, vol.~21, pp. 137-151,
  2005.

\bibitem{Broyden}
{\sc C.~G. Broyden, J.~E.~D. Jr., Wang, and J.~J. More}, {\em On the local and
  superlinear convergence of quasi-newton methods}, IMA J. Appl. Math, 12
  (1973), pp.~223--245.

\bibitem{Bullo2009}
{\sc F.~Bullo, J.~Cortes, and S.~Martinez}, {\em Distributed control of robotic
  networks: a mathematical approach to motion coordination algorithms},
  Princeton University Press, 2009.

\bibitem{byrd2016stochastic}
{\sc R.~H. Byrd, S.~Hansen, J.~Nocedal, and Y.~Singer}, {\em A stochastic
  quasi-{Newton} method for large-scale optimization}, SIAM Journal on
  Optimization, 26 (2016), pp.~1008--1031.

\bibitem{Cao2013-TII}
{\sc Y.~Cao, W.~Yu, W.~Ren, and G.~Chen}, {\em An overview of recent progress
  in the study of distributed multi-agent coordination}, IEEE Transactions on
  Industrial Informatics, 9 (2013), pp.~427--438.

\bibitem{cevher2014convex}
{\sc V.~Cevher, S.~Becker, and M.~Schmidt}, {\em Convex optimization for big
  data: Scalable, randomized, and parallel algorithms for big data analytics},
  IEEE Signal Processing Magazine, 31 (2014), pp.~32--43.

\bibitem{defazio2014saga}
{\sc A.~Defazio, F.~R. Bach, and S.~Lacoste{-}Julien}, {\em {SAGA:} {A} fast
  incremental gradient method with support for non-strongly convex composite
  objectives}, in Advances in Neural Information Processing Systems 27,
  Montreal, Quebec, Canada, 2014, pp.~1646--1654.

\bibitem{gower2016stochastic}
{\sc R.~M. Gower, D.~Goldfarb, and P.~Richt{\'a}rik}, {\em Stochastic block
  {BFGS}: Squeezing more curvature out of data}, arXiv preprint
  arXiv:1603.09649,  (2016).

\bibitem{gurbuzbalaban2015convergence}
{\sc M.~G{\"u}rb{\"u}zbalaban, A.~Ozdaglar, and P.~Parrilo}, {\em On the
  convergence rate of incremental aggregated gradient algorithms}, arXiv
  preprint arXiv:1506.02081,  (2015).

\bibitem{Dennis}
{\sc J.~J.~E.~Dennis and J.~J. More}, {\em A characterization of superlinear
  convergence and its application to quasi-Newton methods}, Mathematics of
  computation, 28 (1974), pp.~549--560.

\bibitem{johnson2013accelerating}
{\sc R.~Johnson and T.~Zhang}, {\em Accelerating stochastic gradient descent
  using predictive variance reduction}, in Advances in Neural Information
  Processing Systems 26, Lake Tahoe, Nevada, United States, 2013, pp.~315--323.

\bibitem{lecun2010mnist}
{\sc Y.~LeCun, C.~Cortes, and C.~J. Burges}, {\em {MNIST} handwritten digit
  database}, AT\&T Labs [Online]. Available: http://yann. lecun.
  com/exdb/mnist,  (2010).

\bibitem{LopesEtal8}
{\sc C.~G. Lopes and A.~H. Sayed}, {\em Diffusion least-mean squares over
  adaptive networks: Formulation and performance analysis}, IEEE Transactions
  on Signal Processing, 56 (2008), pp.~3122--3136.

\bibitem{lucchi2015variance}
{\sc A.~Lucchi, B.~McWilliams, and T.~Hofmann}, {\em A variance reduced
  stochastic {Newton} method}, arXiv preprint arXiv:1503.08316,  (2015).

\bibitem{mokhtari2016surpassing}
{\sc A.~Mokhtari, M.~G{\"u}rb{\"u}zbalaban, and A.~Ribeiro}, {\em Surpassing
  gradient descent provably: A cyclic incremental method with linear
  convergence rate}, arXiv preprint arXiv:1611.00347,  (2016).

\bibitem{mokhtari2014res}
{\sc A.~Mokhtari and A.~Ribeiro}, {\em {RES}: Regularized stochastic {BFGS}
  algorithm}, IEEE Transactions on Signal Processing, 62 (2014),
  pp.~6089--6104.

\bibitem{mokhtari2015global}
{\sc A.~Mokhtari and A.~Ribeiro}, {\em Global convergence of online limited
  memory {BFGS}}, Journal of Machine Learning Research, 16 (2015),
  pp.~3151--3181.

\bibitem{MoritzNJ16}
{\sc P.~Moritz, R.~Nishihara, and M.~I. Jordan}, {\em A linearly-convergent
  stochastic {L-BFGS} algorithm}, in Proceedings of the 19th International
  Conference on Artificial Intelligence and Statistics, {AISTATS} 2016, Cadiz,
  Spain, May 9-11, 2016, 2016, pp.~249--258.

\bibitem{nesterov2013introductory}
{\sc Y.~Nesterov}, {\em Introductory lectures on convex optimization: A basic
  course}, vol.~87, Springer Science \& Business Media, 2013.

\bibitem{Powell}
{\sc M.~J.~D. Powell}, {\em Some global convergence properties of a variable
  metric algorithm for minimization without exact line search}, Academic Press,
  London, UK, 2~ed., 1971.

\bibitem{Ribeiro10}
{\sc A.~Ribeiro}, {\em Ergodic stochastic optimization algorithms for wireless
  communication and networking}, IEEE Transactions on Signal Processing, 58
  (2010), pp.~6369--6386.

\bibitem{Ribeiro12}
{\sc A.~Ribeiro}, {\em Optimal resource allocation in wireless communication
  and networking}, EURASIP Journal on Wireless Communications and Networking,
  2012 (2012), pp.~1--19.

\bibitem{rodomanov2016superlinearly}
{\sc A.~Rodomanov and D.~Kropotov}, {\em A superlinearly-convergent proximal
  newton-type method for the optimization of finite sums}, in Proceedings of
  The 33rd International Conference on Machine Learning, 2016, pp.~2597--2605.

\bibitem{roux2012stochastic}
{\sc N.~L. Roux, M.~W. Schmidt, and F.~R. Bach}, {\em A stochastic gradient
  method with an exponential convergence rate for finite training sets}, in
  Advances in Neural Information Processing Systems 25, Lake Tahoe, Nevada,
  United States., 2012, pp.~2672--2680.

\bibitem{Schizas2008-1}
{\sc I.~Schizas, A.~Ribeiro, and G.~Giannakis}, {\em Consensus in ad hoc wsns
  with noisy links - part i: Distributed estimation of deterministic signals},
  IEEE Transactions on Signal Processing, 56 (2008), pp.~350--364.

\bibitem{Schmidt2016}
{\sc M.~Schmidt, N.~Le~Roux, and F.~Bach}, {\em Minimizing finite sums with the
  stochastic average gradient}, Mathematical Programming,  (2016), pp.~1--30.

\bibitem{SchraudolphYG07}
{\sc N.~N. Schraudolph, J.~Yu, and S.~G{\"{u}}nter}, {\em A stochastic
  quasi-{Newton} method for online convex optimization}, in Proceedings of the
  Eleventh International Conference on Artificial Intelligence and Statistics,
  {AISTATS} 2007, pp.~436--443.

\bibitem{SS}
{\sc S.~Shalev{-}Shwartz and N.~Srebro}, {\em {SVM} optimization: inverse
  dependence on training set size}, in Machine Learning, Proceedings of the
  Twenty-Fifth International Conference {(ICML} 2008), Helsinki, Finland, 2008,
  pp.~928--935.

\end{thebibliography}
\end{document}